\documentclass[a4paper,11pt,english,reqno]{amsart}
\usepackage{mystyle}

\usepackage{bm}
\usepackage{upgreek}
\usepackage{mathabx}

\numberwithin{equation}{section}
\usepackage{url} 

\usepackage{mathrsfs}
\usepackage[normalem]{ulem} 
\usepackage{graphicx,xcolor}

\renewcommand{\Im}{{\rm Im} \,}
\renewcommand{\Re}{{\rm Re} \,}

\title[Large deviations of the empirical spectral measure]{Large deviations of the empirical spectral measure of supercritical sparse Wigner matrices}  
\author{Fanny Augeri} 
\address[Fanny Augeri]{Laboratoire de Probabilit\'es, Statistique et Mod\'elisation (LPSM), Universit\'e Paris Cité, 75205 Paris Cedex 13, France.}
\email{augeri@lpsm.paris}

\begin{document}
\maketitle

\begin{abstract}Let $\Xi$ be the adjacency matrix of an Erd\H{o}s-Rényi graph on $n$ vertices and with parameter $p$ and consider $A$ a $n\times n$ centered random symmetric matrix with bounded i.i.d. entries above the diagonal. When the mean degree $np$ diverges, the empirical spectral measure of the normalized Hadamard product $(A \circ \Xi)/\sqrt{np}$ converges weakly in probability to the semicircle law. In the regime where $p\ll 1$ and $ np \gg  \log n$, we prove a large deviations principle for the empirical spectral measure with speed $n^2p$ and with a good rate function solution of a certain variational problem. The rate function reveals in particular that the only possible deviations  at the exponential scale $n^2p$ are around measures coming from Quadratic Vector Equations. As a byproduct, we obtain a large deviations principle for the empirical spectral measure of supercritical Erd\H{o}s-Rényi graphs. 
 \end{abstract}
 
 \section{Introduction and main result}
\subsection{Setup and main result} Take $A$ to be a  $n\times n$ Rademacher Wigner matrix, meaning that $A$ is a symmetric random matrix such that $(A_{ij})_{i\leq j}$ are independent and uniformly distributed on $\{-1,1\}$. For any $n\times n$ real symmetric matrix $M$, denote by $\mu_M$ its empirical spectral measure, that is $\mu_M:=\frac{1}{n} \sum_{i=1}^n \delta_{\lambda_i}$, where $\lambda_1,\ldots,\lambda_n$ are the eigenvalues of $M$. The celebrated Wigner's Theorem \cite{Wigner} says that in probability, $\mu_{A/\sqrt{n}}$ converges weakly to the semicircle law $\mu_{\text{sc}} := (2\pi)^{-1} \sqrt{(4-x^2)_+}dx$ as $n$ goes to $+\infty$. Now, what is the probability that $\mu_{A/\sqrt{n}}$ is close to a given probability measure different than the semicircle law? By Talagrand's concentration inequality, we know that this probability  decreases exponentially fast with speed $n^2$ \cite{GZconc}. But beyond this fact, no other results - to our knowledge - are known about the large deviation behaviour of $\mu_{A/\sqrt{n}}$. In particular,  the existence of a large deviations principle (LDP) remains open. While many LDPs are known for integrable models of random matrices, like $\beta$-ensembles (\cite{BenArous}, \cite[section 2.6]{AGZ}), $\beta$-Laguerre ensembles \cite{HP}, or general classical ensembles \cite{ES}, the non-integrability of Rademacher Wigner matrices renders the analysis of its large deviation behaviour challenging. Beyond integrable models, the large deviations of the empirical spectral measure are only well-understood in the case of Wigner matrices ``without Gaussian tails'' \cite{BC} by leveraging some heavy-tail phenomenon. Still, impressive results have been obtained recently on the large deviation behaviour of the extreme eigenvalues of Wigner matrices with sub-Gaussian coefficients,  using the novel technique of spherical integrals (see \cite{GH}, \cite{AGH}, \cite{CDG}), which give new hope of solving the problem of the large deviations of the empirical spectral measure of Rademacher Wigner matrices.

Motivated by this question, we consider a sparse relaxation of this model by zeroing out independently each entry on and above the diagonal of $A$ with some probability $p\ll 1$, in such a way that on average, the number of non zero entries on each line and column diverges with $n$, that is, such that $np \to +\infty$. With the appropriate normalization $\sqrt{np}$, the empirical spectral measure still converges weakly to the semicircle law in probability by \cite{KhKhPS}. What can be said of the large deviation behaviour of the empirical spectral measure of this diluted Rademacher matrix? In the present work, we investigate this question more generally  for the model of {\em sparse Wigner matrix with bounded entries} which we describe now. 

\begin{Def} (Sparse Wigner matrix with bounded entries)\label{defmodel}
Let $A$ be a symmetric random matrix of size $n \times n$ with entries $(A_{ij})_{i,j \in [n]}$ (where $[n]:=\{1,2,\ldots, n\}$) such that $(A_{ij})_{i\leq j}$ are i.i.d. bounded random variables with zero mean and unit variance. Let  $p \in (0,1)$ and $\Xi := (\xi_{ij})_{i,j \in [n]}$ be a $n\times n$ symmetric random matrix with zeros on the diagonal and i.i.d. Bernoulli $p$ random variables above the diagonal. We call the matrix $ \widehat{X}:= A  \circ \Xi$, where $\circ$ denotes the Hadamard product between matrices, a {\em sparse Wigner matrix with bounded entries} and set $X:=\widehat{X}/\sqrt{np}$. 

\end{Def}

As we will see, in the regime where $p\ll 1$ and $np \gg \log n$, the mechanisms of large deviations of $\mu_X$  are closely related to changes of the variance profile  of $A\circ \Xi$, and therefore the possible deviations of $\mu_X$ will turn out to be around limiting empirical spectral measures of Wigner matrices with a non-constant variance profile. The measures arising from such Wigner-type matrices are well-understood through their Stieltjes transforms, which can be described as averages of solutions of Quadratic Vector Equations (QVE), studied in depth in \cite{AEK}. 

These QVEs will play a central role in our large deviation analysis and in the variational formula defining our rate function. To present them, we introduce some more notation.
Let $\mathbb{H}$ denote the open upper half plane $\{ z\in \mathbb{C }: \Im z>0\}$, and $\mathcal{P}(\RR)$ the set of Borel probability measures on $\RR$. Denote for any $\mu \in \mathcal{P}(\RR)$ by $m_\mu$ the {\em Stieltjes transform} of $\mu$, defined as 
\begin{equation} \label{defStiel} m_\mu(z) := \int\frac{d\mu(x)}{x-z}, \ z\in \mathbb{H}.\end{equation}
Denote by $\mathcal{B}^+$ the set $\{\ell :[0,1] \to \mathbb{H} : \sup_x|\ell(x)|<+\infty\}$ endowed with the sup norm denoted by $\| \ \|_{\infty}$. Now, consider  a non negative Borel measurable kernel  $W: [0,1]^2\to [0,+\infty)$ which is symmetric in the sense that $W(x,y) =W(y,x)$ for any $(x,y) \in [0,1]^2$, and denote by $d_W$ its {\em degree function}, defined as 
\begin{equation} \label{defdegfunc} d_W(x):=\int_0^1 W(x,y) dy, \ x \in [0,1].\end{equation}
Assuming that $W$ has a bounded degree function,  that is $\sup_{x\in [0,1]} d_W(x) <+\infty$, we know by \cite[Theorem 2.1]{AEK} that for any $z\in \mathbb{H}$, there exists a unique solution $m(z) = (m(z,x))_{x\in [0,1]}$ in $\mathcal{B}^+$ to the Quadratic Vector Equation (QVE) associated to $W$, 
\begin{equation} \label{QVE} - \frac{1}{m(z,x)} = z + \int_0^1 W(x,y) m(z,y) dy, \quad x\in [0,1].\end{equation}
Moreover, the solution $z \mapsto m(z)$ is analytic from $\mathbb{H}$ to $\mathcal{B}_+$ and for each $x$, there exists a probability measure $\upsilon_x \in \mathcal{P}(\RR)$ such that $m_{\upsilon_x}(z) = m(z,x)$ for any $z\in \mathbb{H}$. Further, the map $x\mapsto \upsilon_x$ is Borel measurable. Now, define $\upsilon_W$ as the probability measure
\begin{equation} \label{defup} \upsilon_W:= \int_{[0,1]} \upsilon_x dx.\end{equation}
One can check that if $W'$ is yet another non negative Borel measurable symmetric kernel with a bounded degree function which agrees with $W$ except on a Borel set of zero Lebesgue measure then $\upsilon_W=\upsilon_{W'}$. Setting $\mathcal{X}$ as the set of such kernels, where kernels coinciding almost everywhere for the Lebesgue measure are identified, 
we can thus define unambiguously for any $W\in \mathcal{X}$, the probability measure $\upsilon_W$ by the formula \eqref{defup}. Such measures $\upsilon_W$ where $W\in \mathcal{X}$ arise naturally as the limit of the empirical spectral measures of Wigner-type matrices (see \cite[Theorem 1.1]{Girko}, \cite{Zhu}, \cite{AEK2}). 
To write our rate function, we need to further extend this definition to more general kernels belonging to the set $\mathcal{W}$ described in the following definition. 

\begin{Def}
$\mathcal{W}$ is the set of integrable kernels $W : [0,1]^2\to [0,+\infty)$ which are non negative and symmetric, where kernels agreeing almost everywhere for the Lebesgue measure are identified.
\end{Def}
For 
any $C>0$ and $W\in \mathcal{W}$, denote by $W^{(C)}$ the degree truncated kernel defined by 
\begin{equation} \label{defWC} W^{(C)}(x,y) := W(x,y) \Car_{d_W(x)\leq C} \Car_{d_W(y)\leq C}, \ (x,y) \in [0,1]^2.\end{equation}
Using a generalisation of Hoeffman-Wielandt inequality for measures $\upsilon_{W'}$ arising from kernels $W'\in \mathcal{X}$ (see Lemma \ref{HW}), it follows that for any $W\in \mathcal{W}$, the sequence $(\upsilon_{W^{(C)}})_{C>0}$ converges weakly when $C\to +\infty$, and we define $\upsilon_W$ as the limit. In a nutshell, we can define for any kernel $W\in \mathcal{W}$, a probability measure $\upsilon_W$, which we will call the {\em Quadratic Vector Equation measure of} $W$, or in short the {\em QVE measure} of $W$, as summarised in the following definition.

\begin{Def}[QVE measure of a kernel]\label{specW}
For any $W\in \mathcal{X}$, the {\em QVE measure} of $W$ is the probability measure $\upsilon_W :=\int_0^1 \upsilon_x dx$, where $x\mapsto \upsilon_x$ is a Borel measurable function from $[0,1]$ to $\mathcal{P}(\RR)$ endowed with the weak topology such that for any $z\in \mathbb{H}$, $m(z)=(m(z,x))_{x\in [0,1]}$, where $m(.,x)$ is the Stieltjes transform of $\upsilon_x$, is the unique solution in $\mathcal{B}^+$ of the QVE \eqref{QVE} associated to $W$.
This definition extends to any kernels $W\in \mathcal{W}$ by setting $\upsilon_W:=\lim_{C\to+\infty} \upsilon_{W^{(C)}}$ where $W^{(C)}$ is defined in \eqref{defWC}.

\end{Def}

We introduce further the  functions $L$ and $h_L$  defined as,
\begin{equation} \label{defL} L(\theta) = \EE(e^{\theta A_{12}^2})-1, \theta \in \RR,\ \text{ and } \ h_L(u) = \sup_{\theta \in \RR}\{\theta u - L(\theta)\}, u\in \RR.\end{equation}
With this notation, we can now define the functional $I_L$ on $\mathcal{P}(\RR)$ by
\begin{equation} \label{defJ} I_L(\mu) :=  \inf\big\{ \frac{1}{2}\int_{[0,1]^2} h_L(W(x,y)) dx dy : W \in \mathcal{W},  \upsilon_W =\mu\big\},\quad \mu \in \mathcal{P}(\RR),\end{equation}
where $h_L$ is defined in \eqref{defL} and $\upsilon_W$ is defined in Definition \ref{specW} for any $W\in \mathcal{W}$. Our main result reads as follow.
\begin{The}\label{main}
\label{main}Let $\widehat{X}$ be a sparse Wigner matrix with bounded entries and set $X:=\widehat{X}/\sqrt{np}$. Assume $p\ll 1$ and $np \gg \log n$.
The sequence $\mu_X$ satisfies a large deviation principle in $\mathcal{P}(\RR)$ endowed with the weak topology, with speed $n^2p$ and good rate function $I_L$.
Moreover, for any $\mu \in \mathcal{P}(\RR)$ such that $I_L(\mu)<+\infty$ the infimum defining $I_L(\mu)$ is achieved.

\end{The}

A close inspection of the proof reveals that the same arguments as for the empirical spectral measure of $A \circ \Xi /\sqrt{np}$, where the entries of $A$ are Rademacher distributed, hold as well for the large deviations of the empirical spectral measure of $(\Xi-\EE(\Xi))/\sqrt{np}$ with the same rate function.  Using Cauchy interlacing inequality (see \cite[Theorem 1.43]{BS}), this yields the following result on the large deviations of the empirical spectral measure of sparse Erd\H{o}s-Rényi graphs. 
\newpage
\begin{The}\label{LDPER}Let $\Xi$ be the adjacency matrix of an Erd\H{o}s-Rényi graphs on $n$ vertices and with parameter $p$ such that $p\ll 1$ and $np \gg \log n$. 
The sequence $(\mu_{\Xi/\sqrt{np}})_{n\in\NN}$ satisfies a LDP for the weak topology with speed $n^2p$ and with good rate function $I$ defined by, 
\[ I(\mu) :=  \inf\big\{ \frac{1}{2}\int_{[0,1]^2} h(W(x,y)) dx dy : W \in \mathcal{W},  \upsilon_W =\mu\big\},\quad \mu \in \mathcal{P}(\RR),\]
where $h(u) = u \log u-u+1$ for any $u\geq 0$. 
\end{The}

To complement Theorems \ref{main} and \ref{LDPER}, we check that the semicircle law, as typical limit of the empirical spectral measure of supercritical sparse Wigner, is the unique minimizer of the rate function $I_L$.

\begin{Lem}\label{zeroI} Let $\mu \in \mathcal{P}(\RR)$.
$I_L(\mu) = 0$ if and only if $\mu$ is the semicircle law.
\end{Lem}

Finally, we close this section on several remarks on the rate function we obtained and a discussion of our assumptions. 
\begin{Rem}
\begin{enumerate}
\item 
The form of the rate function comes from the fact that the optimal changes of measures are $\propto \exp\big(\sum_{i<j} \theta_{ij} \xi_{ij} A_{ij}^2\big) d\PP$, where $(\theta_{ij})_{i,j}$ is some symmetric matrix corresponding to a stepped admissible kernel. More precisely,  $\theta_{ij} = n^2 \int_{I_{in}\times I_{jn}} \theta(x,y) dxdy$, $i,j \in [n]$, where $\theta : [0,1]^2 \to \RR$ is some admissible kernel and $I_{in}:= (\frac{i-1}{n},\frac{i}{n}]$, $i \in [n]$. Under such a change of measure, $A\circ \Xi$ is centered, has independent entries above the diagonal, and one can check that its variance profile is approximately $(pL'(\theta_{ij}))_{i,j\in[n]}$. By \cite[Theorem 1.1]{Girko}, the resulting deviation of the empirical spectral measure of $X$  will be around the probability measure $\upsilon_W$ where  $W= L' \circ \theta$. Specifying this strategy in the case where $A$ is a Rademacher Wigner matrix, this means that to create a deviation of the empirical spectral measure of $X$ around $\upsilon_W$ where $W$ is well-behaved (for example bounded), it suffices to change the distribution of $\Xi$ into the one of a symmetric matrix with independent coefficients above the diagonal, such that  $\xi_{ij}$ has a Bernoulli $p W_{ij}$ distribution, where $W_{ij} =n^2\int_{I_{in}\times I_{jn}} W(x,y) dx dy$, for any $i,j\in[n]$. As remarked before the statement of Corollary \ref{LDPER}, the same remains true for the large deviations of the empirical spectral measure of supercritical Erd\H{o}s-Rényi graphs. This means that the optimal large deviation strategy is to change the distribution of the graph into one of an inhomogenous Erd\H{o}s-Rényi graph associated to a  well-behaved kernel. 
\item An interesting consequence of Theorem \ref{main} is that the only possible deviations of the empirical spectral measure of sparse Wigner matrices with bounded entries are around symmetric probability measures $\mu$, that is, such that $\mu(-E)=\mu(E)$ for any Borel subset $E$ of $\RR$. Indeed,  from Theorem \ref{main}, we know that if $I_L(\mu)<+\infty$ then $\mu$ should be the QVE measure of some kernel $W\in \mathcal{W}$. Now, by \cite[Theorem 2.1]{AEK} we know that for any $W'\in \mathcal{X}$, the probability measures $\upsilon_x$, $x \in [0,1]$, coming from the QVE measure associated to $W'$, are symmetric which implies that  $\upsilon_{W'}$ is also symmetric. Clearly, this property immediately extends to kernels in $\mathcal{W}$, by definition of the QVE measure (see Definition \ref{specW}), and as a result $\mu$ has to be a symmetric measure. 
\item Our large deviations principle does not apply in the entire regime where the empirical spectral measure converges typically to the semicircle law, that is $np \gg 1$, but only in the supercritical regime where $np \gg \log n$. While the speed of large deviations should remain the same, we do not believe the rate function $I_L$ to be valid when $1\ll np \lesssim \log n$.  Our strategy of proof breaks down in two major points, which are highlighted in the outline of proof (see section \ref{outline}) for the interested reader.
\item For sub-Gaussian, but unbounded entry distributions, in addition to the strategy of changing the variance profile, a new scenario can emerge in the large deviation behaviour of the empirical spectral measure: one can have $O(n)$ entries of order $1$, that is $\xi_{ij}=1$ and $A_{ij}\asymp \sqrt{np}$, which would create a deviation of the empirical spectral measure around a free convolution with the semicircle law, and which has a cost at the exponential scale of order $n^2p$. We believe that this ``heavy-tail scenario'', already found to happen in the case of Wigner matrices ``without Gaussian tails'' \cite{BC}, should coexist with the one of changing the variance profile, in the large deviations of the empirical spectral measure of sparse sub-Gaussian Wigner matrices.
\end{enumerate}

\end{Rem}

\subsection{Related work}Our approach to the large deviations of the empirical spectral measure of sparse Wigner matrices borrows very much from the random graphs literature. 
Large deviations of sparse random graphs have attracted in the past decades considerable attention. Much of the efforts were devoted to understand the large deviations of subgraph counts or homomorphisms densities in sparse Erd\H{o}s-Rényi graphs. Unlike the dense case \cite{CV}, the lack of ``objective method'' or limit theory for sparse graphs renders the study of the large deviations of these observables challenging, the {\em infamous upper tail problem} of triangle counts \cite{JR} being once the epitome of this difficulty imposed by sparsity. After a long sequence of works settling first the speed of deviation (see  \cite{KV}, \cite{JOR}, \cite{Ctriangle}, \cite{DK}), the large deviation rate function of triangle counts was then identified as the solution of a certain mean-field variational problem in \cite{CD16} and \cite{LZ} for sparsity parameters decreasing polynomially fast with the number of vertices. The range of sparsity was successively improved by \cite{Eldan}, \cite{CoD}, \cite{NLAu} until the large deviations in the whole supercritical regime were finally settled in \cite{HMS22}. From these works, roughly three different approaches emerged, the first being the so-called ``nonlinear large deviation theory'' initiated by Chatterjee and Dembo \cite{CD16} and refined by \cite{Eldan}, \cite{A}, \cite{NLAu}, \cite{transpA}, which has the advantage to be very general but usually yields suboptimal results in terms of sparsity range. The second approach \cite{CoD}, closer in spirit to the original  work of Chatterjee and Varadhan \cite{CV} for dense Erd\H{o}s-Rényi graphs, revisits the Regularity Method in the sparse setting and was developed further in \cite{CoDP} to investigate the large deviations of homomorphism densities in Erd\H{o}s-Rényi hypergraphs. Finally, the third approach is a combinatorial method introduced in \cite{HMS22} which has proven to be powerful for understanding localisation phenomena in the large deviations of random graphs, and  was further pursued in \cite{BB} where the authors computed the large deviation upper tail of connected regular subgraph counts of sparse Erd\H{o}s-Rényi graphs.

Regarding large deviations of the spectrum of sparse random graphs, most of the available results study the atypical behaviour of the extreme eigenvalues. For Erd\H{o}s-Rényi graphs $G(n,p)$, the top eigenvalue is typically asymptotically equivalent when $n$ goes to infinity to the maximum of the mean degree and the square root of the maximum degree by \cite{KS}. The transition between these two behaviours occurs when $np \asymp \sqrt{\log n/\log \log n}$. In the regime where the top eigenvalue is typically equivalent to the mean degree, the large deviation upper tail has been investigated in \cite{CoD} and \cite{B}. Yet in this sparsity regime, another transition happens when $np \asymp \log n$ in the typical behaviour of the second largest eigenvalue or of the top eigenvalue of the recentered adjacency matrix (see \cite{ADK21, TY21}). When $np \gg \log n$, these spectral observables stick to the right edge of support of the semicircle law and their large deviation upper tails have been computed in \cite{AB} (although with the additional restriction that $\log (np) \gtrsim \log n$ for the second largest eigenvalue). Regarding the sparsity range $\sqrt{\log n/\log \log n } \ll np \lesssim \log n$, the large deviation behaviour of the second largest eigenvalue or of the top eigenvalue of the recentered adjacency matrix remains open.

In the complement regime where $np \ll \sqrt{\log n/\log \log n}$ and the top eigenvalue of $G(n,p)$ is typically equivalent to the square root of the maximum degree, the joint large deviation upper tails of the extreme eigenvalues were obtained in \cite{BBG}. Moreover, for the top eigenvalue of Erd\H{o}s-Rényi graphs with constant mean degree and Gaussian conductances, a precise description of the upper and lower large deviations tails was derived in \cite{GN}, and large deviation tails for more general distributions of conductances were computed in \cite{GHN22}.

In contrast, much less is known of the large deviations of global observables of the spectrum like the empirical spectral measure. The only known result to our knowledge is for Erd\H{o}s-Rényi graphs with constant mean degree. Indeed, a large deviations principle was proven for this model by Bordenave and Caputo \cite{BC} with respect to the local weak topology and with a good rate function. As the expected spectral measure is continuous for this topology (see \cite[Lemma 3.13]{BC}, \cite[Theorem 4]{ATV}), the contraction principle implies that the empirical spectral measure of Erd\H{o}s-Rényi graphs with constant mean degree satisfies indeed a large deviations principle with respect to the weak topology.

\subsection{Outline of the proof} \label{outline}
From the different approaches to the large deviations of sparse random graphs mentioned in the previous section, one can say that on a high level ours bears most similarities with the Regularity method, as put forward in \cite{CV}, \cite{CoD}, \cite{CoDP}.  Our strategy is to ultimately contract a LDP for a certain empirical kernel associated to our random matrix $X$ with respect to the cut metric. Denote by $G$ the resolvent of $X$ defined by $G(z):= (X-z)^{-1}$ for any $z\in \mathbb{H}$. In the sequel we will often drop the $z$-dependence to make the notation more concise. The starting point of our analysis is the Schur complement formula \cite[(4.1)]{BGK}, which states that for any $z\in \mathbb{H}$,
\begin{equation} \label{Schurform}-\frac{1}{G_{ii}} = z +\sum_{k,\ell}^{(i)} X_{ik}X_{i\ell} G_{k\ell}^{(i)}, \ i\in [n],\end{equation}
where $G^{(i)}$ is the resolvent of the matrix $X^{(i)}$ obtained from $X$ by removing the $i^{\text{th}}$ line and column, and $\sum_{k,\ell}^{(i)}$ denotes the sum $\sum_{k,\ell \neq i}$. One can rewrite \eqref{Schurform} as a perturbed QVE by singling out the diagonal terms on the right-hand side as follows
\begin{equation} \label{pertQVE} -\frac{1}{G_{ii}} = z + \sum_{k=1}^n X_{ik}^2 G_{kk} + d_i(z), \ i\in [n],\end{equation}
where $d_i(z):= \sum_{k\neq i} X_{ik}^2(G_{kk}^{(i)}-G_{kk}) + \sum_{k\neq \ell} X_{ik}X_{i\ell} G_{k\ell}^{(i)}$ for any $i\in[n]$.  On the one hand, using resolvent identities and the fact that the entries of $A$ are bounded, the first term in the perturbation $d_i(z)$ can be shown to be of order $O(1/np)$.  On the other hand, the second term $\sum_{k,\ell}^{(i)} X_{ik}X_{i\ell} G_{k\ell}^{(i)}$ can be seen as a chaos of order 2 in the independent sparse random variables $(X_{ik})_{k\in [n]}$ by conditioning on $X^{(i)}$. Since the sum involves the products $\xi_{ik}\xi_{i\ell}$, $k,\ell \neq i$, which are Bernoulli variables of parameter $p^2$ and we assumed $p\ll 1$, one can show (see Lemma \ref{chaos}) that the probability that such a sum is of order $1$ is exponentially small with a speed much larger than $np$. If the variables $(d_i(z))_{i\in [n]}$ were independent, we would immediately have using Bennett's inequality \cite[Theorem 2.9]{BLM} and the fact that $np \gg \log n$, that at the exponential scale $n^2p$ only an arbitrarily small proportion of $d_i(z)$'s can be non-trivial. To circumvent the issue posed by the non-independence, we propose a generalisation of Bennett's inequality to dependent variables, similar to the one proven by Chatterjee \cite[Theorem 3.1]{Chtriangles} in the context of finding the order of the large deviation upper tail of triangle counts in sparse Erd\H{o}s-Rényi graphs. As a result, we obtain in Proposition \ref{negli} that only a small fraction of the variables $d_i(z)$, $i\in[n]$ can be non-trivial at the large deviation scale $n^2p$.
This result prompts us to introduce the random edgeweighted graph $\mathcal{G}_n$ with adjacency matrix $(A_{ij}^2\xi_{ij})_{i,j\in [n]}$ and define $W_n$ the associated kernel on $[0,1]^2$ by setting $W_n$ to be equal to $A_{ij}^2\xi_{ij}/p$ on the square $(\frac{i-1}{n},\frac{i}{n}]\times (\frac{j-1}{n},\frac{j}{n}]$ for any $i,j \in [n]$.
We can then veritably view \eqref{pertQVE} at the large deviation scale $n^2p$ as a small perturbation of the following QVE 
\begin{equation} \label{unpertQVE} -\frac{1}{m_i} = z + \sum_{k=1}^n X_{ik}^2 m_k, \ i\in[n]\end{equation}
which is the QVE associated to the empirical kernel $W_n$: as $W_n$ is a stepped kernel, the corresponding QVE \eqref{QVE} simplifies itself into the finite-dimensional equation \eqref{unpertQVE} on $\mathbb{H}^n$. Now, building on the stability results of QVE obtained by Ajanki, Erd\H{o}s and Kruger  in \cite{AEK}, we deduce in Lemma \ref{equivexpo}
 that $(\upsilon_{W_n})_{n\in \NN}$ is an exponential equivalent of $\mu_X$, in the sense of \cite[Definition 4.2.10]{DZ}.

At this point, we are in a favourable position to apply the contraction principle (see \cite[Theorem 4.2.1]{DZ}). First, we embed the empirical kernels $(W_n)_{n\in\NN}$ into the space of  integrable non negative symmetric kernels $\mathcal{W}$, which we equip of the topology induced by the cut norm. An immediate difficulty arising from such an embedding is to control the complexity of the kernel $W_n$ for the cut norm, so that one can reduce  the large deviation upper bound to computing ball probabilities. To this end, we show in Proposition \ref{hypoH} that the kernel $W_n$ is upper regular in the terminology of \cite{BCCZ} with overwhelming probability. This entails using the Regularity Lemma of Borgs, Chayes, Cohn and Zhao \cite[Theorem C.11]{BCCZ} proven for upper regular kernels that with overwhelming probability the kernel $W_n$ lives in a subset of $\mathcal{W}$ of metric entropy relative to the cut norm at most $n \log n$ (see Proposition \ref{expopresqtens}). Since $np \gg \log n$, the complexity $W_n$ for the cut norm is thus negligible compared to our large deviation speed $n^2p$. Equipped with such a result, we derive in Proposition \ref{LDPW} a LDP for $(W_n)_{n\in \NN}$ with respect to the cut norm. We then contract this LDP to the space of {\em unlabelled kernels} $\widetilde{\mathcal{W}}$ and obtain a LDP for the unlabelled kernels $(\widetilde{W}_n)_{n\in \NN}$. Making use of the compactness for the cut metric of uniformly integrable kernels proven in \cite[Theorem C.7]{BCCZ}, we show that the resulting rate function of the LDP of $(\widetilde{W}_n)_{n\in \NN}$ has compact level sets, a key assumption in the contraction principle.

 We now move on to check the continuity of our observable - the QVE measure of a kernel - with respect to the cut norm. When the kernel has a bounded degree function, its spectral measure is a symmetric, compactly supported probability measure by \cite[Theorem 2.1]{AEK} and its even moments are expressed in terms of homomorphism densities of trees by \cite[Lemma 2.2]{EM} (see Proposition \ref{momentuW}). Thus, we prove a Counting Lemma \ref{countinglem} restricted to homomorphism densities of trees for kernels with bounded degree function, and use this result to show in Proposition \ref{conti} that the map $W\mapsto \upsilon_W$ is uniformly continuous on the set $\mathcal{X}_C$ of kernels with degree function uniformly bounded by $C$, for any given $C>0$. In Proposition \ref{defcont} we are able to leverage this continuity to the whole space $\mathcal{W}$. Finally, as the QVE measure of a kernel is invariant by relabelling and continuous for the cut norm, we show that we can define the QVE measure of any unlabelled kernel and that this map is continuous on $\widetilde{\mathcal{W}}$ for the cut metric (see Proposition \ref{conttilde}). All the requirements of the contraction principle are now met, which ends the sketch of the proof of Theorem \ref{main}.

\subsection{Notation}
For any set $T$, we denote by $\# T$ its cardinal. For any vector $u\in \RR^n$, we denote by $\|u\|$ its $\ell^2$ norm. For any $n\times n$ matrix $M$ we denote by $\|M\|$ its operator norm with respect to the $\ell^2$ norm on $\RR^n$, by ${\rm Tr(M)}$ its trace and ${\rm rank}(M)$ its rank. Further, we denote by $M_j$ its $j^{\text{the}}$ column for any $j\in [n]$, and for any  $T\subset [n]$,  by $M^{(T)}$ the $T^c\times T^c$ submatrix of  $M$ spanned by the lines and columns in $T^c$. For any rectangular matrix $M$, we denote by $M^{\sf T}$ its transpose matrix. For any functions $\theta,W :[0,1]^2\to\RR$ such that $\theta W\in L^1([0,1]^2)$, we denote by $\langle \theta, W\rangle = \int_{[0,1]^2} \theta(x,y) W(x,y) dx dy$. Finally, we denote by $\lambda$ the Lebesgue measure on $[0,1]$ and $\lambda^2$ the Lebesgue measure on $[0,1]^2$. 

\section{An exponential equivalent}
In the rest of this paper, we denote by $R$ the essential supremum of the $A_{ij}$'s for any $i,j \in [n]$. Our first major step is to prove that $\mu_X$ is exponentially equivalent to the QVE measure of a certain empirical kernel. 
First, with the sparsity parameter $p$ being fixed, we associate to any edgeweighted graph a kernel as follows. 

\begin{Def}\label{WG}
For any edgeweighted graph $\mathcal{G}$ with vertex set $[n]$ and adjacency matrix $\beta$, the {\em kernel associated to $\mathcal{G}$}, denoted by $W^\mathcal{G}$, is defined as the function in $\mathcal{W}$ taking the value $\beta_{ij}/p$ on each square $I_{in}\times I_{jn}$ for any $i,j \in [n]$, where $I_{k n} = (\frac{k-1}{n},\frac{k}{n}]$ for any $k\in[n]$.  
\end{Def}
Note that the definition of the kernel $W^\mathcal{G}$ depends implicitly on the sparsity parameter $p$. This normalisation comes from the fact that all the graphs we will consider in this paper have an edge density of order $p$. Set 
$\mathcal{G}_n$ to be the random edgeweighted graph with adjacency matrix $\beta_n:=(A_{ij}^2\xi_{ij})_{i,j\in [n]}$ and denote by  $W_n:=W^{\mathcal{G}_n}$  the associated kernel in $\mathcal{W}$. More precisely, for any $(x,y) \in [0,1]^2$ we have
\begin{equation} \label{defWn} W_n(x,y) := \frac{\xi_{ij}}{p} A_{ij}^2, \quad   (x,y) \in I_{in}\times I_{jn} .\end{equation}
The statement and proof of our exponential equivalent will involve different metrics on $\mathcal{P}(\RR)$ that we introduce now. Since we will have to handle Stieltjes transforms, it will be convenient to work with  the following distance $d$ on $\mathcal{P}(\RR)$, compatible with the weak topology,
 \begin{equation} \label{defd} d(\mu,\nu) := \sup \big\{|m_\mu(z) -m_\nu(z)| : \Im z \geq 2, z \in \mathbb{H}\big\}, \  \mu,\nu \in \mathcal{P}(\RR),\end{equation}
 where $m_\mu$, $m_\nu$ are the Stieltjes transforms of $\mu$ and $\nu$ respectively, defined in \eqref{defStiel}.  Further, let $d_{\text{KS}}$ denote the {\em Kolmogorov-Smirnov distance} defined as
 \begin{equation} \label{defKS} d_{\text{KS}}(\mu,\nu) := \sup_{t\in \RR} | F_\mu(t) - F_\nu(t)|, \ \mu,\nu \in \mathcal{P}(\RR)\end{equation}
 where $F_\mu(t):=\mu((-\infty,t])$, $t\in \RR$, denotes the distribution function of $\mu$, and similarly for $F_\nu$. Finally, let $\mathscr{W}_p$,  $p\geq 1$, denote the {\em $L^p$ Wasserstein distance}, 
  \begin{equation} \label{defWp}   \mathscr{W}_p(\mu,\nu) :=\inf_\pi \Big(\int |x-y|^p d\pi(x,y)\Big)^{1/p}, \ \mu,\nu\in \mathcal{P}(\RR),\end{equation}
 where the infimum runs over couplings $\pi$ between $\mu$ and $\nu$. 
The metric $d$ is related to the $L^1$ Wasserstein distance and the Kolmogorov-Smirnov distance by the inequality (see \cite[(14)]{BC}):
\begin{equation} \label{ineqd} d(\mu,\nu) \leq \mathscr{W}_1(\mu,\nu) \wedge d_{\text{KS}}(\mu,\nu), \ \mu,\nu \in \mathcal{P}(\RR),\end{equation}
which we will use often in the sequel. With this notation, the goal of this section is to prove the following exponential equivalent.


\begin{Pro}\label{equivexpomain1}

For any $\veps>0$,
\[ \lim_{n\to+\infty} \frac{1}{n^2p} \log \PP\big( d(\mu_X,\upsilon_{W_n})>\veps\big) =-\infty,\]
where $d$ is defined in \eqref{defd}, and $\upsilon_{W_n}$ in Definition \ref{specW}.
\end{Pro}

Note that as $W_n$ is a stepped kernel, the definition of $\upsilon_{W_n}$ simplifies itself in the following way.  We know by \cite[Theorem 2.1]{AEK} that for any $z\in \mathbb{H}$ there is a unique solution $m(z):=(m_i(z))_{i\in [n]}$ in $\mathbb{H}^n$ to the QVE
\begin{equation}  - \frac{1}{m_i(z)} = z+\sum_{i=1}^n X_{ij}^2 m_j(z), \ i \in[n]\end{equation}
and that for any $i\in [n]$, $m_i$ is the Stieltjes transform of a probability measure $\upsilon_i$ on $\RR$.
Using the unicity of the solution of the QVE \eqref{QVE} associated to the stepped kernel $W_n$ garanteed again by \cite[Theorem 2.1]{AEK}, it follows that the solution for any $z\in \mathbb{H}$ is piecewise constant and given by ${m}(z,x) =m_{\lceil nx\rceil}(z)$ for any $x\in [0,1]$. As a consequence, we have the simpler representation of the QVE measure of $W_n$ as  $\upsilon_{W_n} = \frac{1}{n} \sum_{i=1}^n \upsilon_i$. This fact holds for the QVE measure of any stepped kernel, which we will use repetitively in this section. 

To prove Proposition \ref{equivexpo}, we will need a few intermediate results. 
In a first step, we show that it is sufficient to consider the submatrix of $X$ spanned by lines and columns with bounded $\ell^2$ norms, for which it will be technically easier to prove the exponential equivalence. 
To prove this reduction, we will show that the proportion of columns (or lines) with large $\ell^2$ norm is negligible at the exponential scale $n^2p$.  More precisely, denote for any $i\in [n]$ by $X_i$ the $i^{\text{th}}$ column of $X$. Now, for any $C\geq 1$
let 
\begin{equation} \label{defJC}{J}_C:=\{ i\in [n] : \|X_i\|\geq C\}.\end{equation}
In the following lemma, we show  that at the exponential scale $n^2p$, ${J}_C$ carries only a small fraction of the total number of columns as $C\to+\infty$. 

\begin{Lem}\label{deg}
For any $t>0$, 
\[\lim_{C\to+\infty} \limsup_{n\to+\infty} \frac{1}{n^2p} \log \PP\big( |J_C|\geq t n\big) = -\infty.\]
\end{Lem}

\begin{proof}Let $C>1$. 
Define for any $i\in [n]$, $Z_i:= {\bf 1}{\{\sum_{j>i} \xi_{ij} \geq C^2np\}}$ and set $v:= \sum_{i=1}^n \EE( Z_i)$. As the random variables $\{(\xi_{ij})_{j>i}\}_{i\in [n]}$ are independent, by Bennett's inequality \cite[Theorem 2.9]{BLM} we know that for any $t n >v$, 
\begin{equation} \label{Ben1} \PP\big( \sum_{i=1}^n Z_i> t n\big) \leq \exp\Big(- t n \log  \frac{t n}{3v} \Big),\end{equation}
Now, using Bennett's inequality again gives that for any $i\in [n]$ 
\begin{equation}\label{taildeg} \EE( Z_i) \leq  \PP\big( \sum_{i=1}^n \xi
_{ij}\geq C^2np \big) \leq \exp\big(-np C^2 \log (C^2/3)\big).\end{equation}
As $np \gg \log n$, for $n$ large enough $v = \sum_{i=1}^n \EE (Z_i) \leq \exp(-np (C^2/2)\log(C^2/3))$. Coming back to \eqref{Ben1}, this gives for any $t>0$ and $n$ large enough, 
\begin{equation} \label{boundbinom} \PP\big( \sum_{i=1}^n Z_i > t n \big) \leq \exp\Big( - t n^2 p (C^2/3) \log (C^2/3)\Big).\end{equation}
Symmetrically, the same bound holds for $\sum_{i=1}^n \widetilde{Z}_i$, where $\widetilde{Z}_i = {\bf 1}{\{\sum_{i<j} \xi_{ij} \geq C^2np\}}$ for any $i \in[n]$.  As the entries of $A$ are bounded by $R$, we have 
\[\#\{i:\|X_i\|^2\geq  2C^2R^2\} \leq \#\{ i : \sum_{i=1}^n \xi_{ij}\geq 2C^2np\} \leq \sum_{i=1}^n Z_i + \sum_{i=1}^n \widetilde{Z}_i.\]
Using \eqref{boundbinom} and a union bound, this gives the claim. 
\end{proof}

Denote for any $T \subset [n]$ by $X^{(T)}$ the $T^c \times T^c$ submatrix of $X$ spanned by the lines and columns in $T^c$. Recall Cauchy interlacing inequality (see  \cite[Theorem 1.43]{BS}) which states that for any two $n\times n$ real symmetric matrices $M,M'$, 
\begin{equation} \label{rkineqm} d_{KS}(\mu_M,\mu_{M'}) \leq \frac{1}{n} {\rm rank}(M-M').\end{equation}
Using Lemma \ref{deg} and the above inequality, 
it follows immediately that the empirical spectral measure of the submatrix $X^{(J_C)}$ is an exponentially good approximation of $\mu_X$, in the sense of \cite[Definition 4.2.14]{DZ}. More precisely, we have the following lemma. 

\begin{Lem} \label{reduc}For any $\veps>0$,
\[ \lim_{C\to+\infty} \limsup_{n\to+\infty} \frac{1}{n^2p} \log \PP\big(d_{\text{KS}}(\mu_X,\mu_{X^{(J_C)}})>\veps\big) = -\infty,\]
where $J_C$ is defined in \eqref{defJC}.

\end{Lem}
\begin{proof}To ease the notation, we write $J$ instead of $J_C$. 
Let $\widetilde{X}^{(J)}$ be the matrix obtained from $X$ by zeroing out the lines and columns in $J$. Since the rank of $X-\widetilde{X}^{(J)}$ is at most $2|J|$,  we deduce by the rank inequality \eqref{rkineqm} that $d_{\text{KS}}(\mu_X,\mu_{\widetilde{X}^{(J)}})\leq 2|J|/n$. Besides, $\mu_{\widetilde{X}^{(J)}}=(|J^c|/n) \mu_{X^{(J)}} + (|J|/n) \delta_0$, so that $d_{\text{KS}}(\mu_{X^{(J)}},\mu_{\widetilde{X}^{(J)}})\leq 2|J|/n$. By the triangle inequality it follows that $d_{\text{KS}}(\mu_X,\mu_{{X}^{(J)}}) \leq 4|J|/n$. Invoking Lemma \ref{deg}, this ends the proof of the claim.

\end{proof}
Similarly, let $W_n^{(J_C)}$ be the kernel associated to the edgeweighted graph with vertex set $J_C^c$ and adjacency matrix $\beta_n^{(J_C)}$. 
As in Lemma \ref{reduc}, we show that $\upsilon_{W_n^{(J_C)}}$ is an exponentially good approximation of $\upsilon_{W_n}$. 
\begin{Lem}\label{approxW}
For any $\veps>0$, 
\[ \lim_{C\to+\infty} \limsup_{n\to+\infty} \frac{1}{n^2p} \log \PP\big(d(\upsilon_{W_n},\upsilon_{W_n^{(J_C)}}) >\veps\big) =-\infty,\]
where $d$ is defined in \eqref{defd}.
\end{Lem}

To prove Lemma \ref{approxW}, we will use the following  version of Cauchy interlacing inequality for QVE measures of kernels.

 \begin{Lem}[Cauchy interlacing inequality]\label{rkineq}
Let $W,W'\in \mathcal{X}$. Assume that there exists a Borel subset $E$ of $[0,1]$ such that ${W'}(x,y) = W(x,y)$ for almost all $(x,y)$ in $E^c\times E^c$. Then, 
\[ d(\upsilon_W,\upsilon_{W'}) \leq 2\lambda(E),\]
where $d$ is defined in \eqref{defd}, and $\lambda$ denotes the Lebesgue measure on $[0,1]$.
\end{Lem}

The proof of this inequality can be found in the Appendix \ref{appendix}. Equipped with this result, we are now ready to give a proof of Lemma \ref{approxW}.

\begin{proof}[Proof of Lemma \ref{approxW}] Writing again $J$ instead of $J_C$, let $\widetilde{\beta}_n^{(J)}$ be the matrix obtained from $\beta_n$ by zeroing out the lines and columns in $J$, and denote by $\widetilde{W}_n^{(J)}$ the associated kernel. Let $g$, $\widetilde{g}$ denote the Stieltjes transforms of $\upsilon_{W_n^{(J)}}$ and $\upsilon_{\widetilde{W}_n^{(J)}}$ respectively. By definition,  $g = \sum_{i \in J^c} m_i /|J^c|$, and $\widetilde{g}= \sum_{i=1}^n \widetilde{m}_i/n$, where for any $z\in \mathbb{H}$, $m(z): =(m_i(z))_{i\in J^c}$ and $\widetilde{m}(z):=(m_i(z))_{i\in[n]}$ are the unique solutions in $\mathbb{H}^{J^c}$ and $\mathbb{H}^n$ respectively of the QVEs associated to $W_n^{(J)}$ and $\widetilde{W}_n^{(J)}$ respectively, that is, 
\begin{equation} \label{QVEX} -\frac{1}{m_i(z)} = z + \sum_{j}^{(J)} X_{ij}^2 m_j(z), \ i\in J^c,\end{equation}
and 
\begin{equation} \label{QVEX2} -\frac{1}{\widetilde{m}_i(z)} = z + \sum_{j}^{(J)} X_{ij}^2 \widetilde{m}_j(z), \ i\in J^c, \quad \widetilde{m}_i(z)= -\frac{1}{z}, \ i \in J.\end{equation}
 By unicity of the solution of the equation \eqref{QVEX}, it follows that for any $z\in \mathbb{H}$ and $i\in J^c$, $\widetilde{m}_i(z) = m_i(z)$. Thus, for any $z\in \mathbb{C}$, $\Im z \geq 2$, 
\[ |g(z)- \widetilde{g}(z)| \leq \frac{|J|}{n} |g(z)| + \frac{|J|}{n|z|} \leq \frac{|J|}{n},\]
where we used the fact that $|g(z)|\leq 1/\Im z$. As a consequence $d(\upsilon_{W_n^{(J)}},\upsilon_{\widetilde{W}_n^{(J)}}) \leq |J|/n$. Now, using Lemma \ref{rkineq}, we know that $d(\upsilon_{W_n},\upsilon_{\widetilde{W}_n^{(J)}}) \leq 2|J|/n$. By the triangular inequality, this yields $d(\upsilon_{W_n},\upsilon_{W_n^{(J)}})\leq 3|J|/n$. Using Lemma \ref{deg}, this ends the proof.

\end{proof}

Putting together Lemmas \ref{reduc} and \ref{approxW}, it follows that in order to prove Proposition \ref{equivexpomain1}, it is sufficient to show, for $C$ large enough, that $\mu_{X^{(J_C)}}$ and $\upsilon_{W_n^{(J_C)}}$ are arbitrary close with overwhelming probability. 

\subsection{A concentration argument}
To ease the notation,  we fix in the following a constant $C\geq 1$ large enough and drop the dependency of the set $J_C$ in $C$, and write $J$ instead of $J_C$ to denote the set of columns with $\ell^2$ norm greater than $C$. We will prove that the Stietljes transform of $\mu_{X^{(J)}}$ is the average of the solution of a perturbation of the  QVE associated to $W_n^{(J)}$. 
Define $G$ the resolvent of $X$ by
\[ G(z) := (X-z)^{-1}, \ z \in \mathbb{H}.\]
When there is no risk of confusion, we will drop the dependency of the resolvent in $z$ and write $G$ instead of $G(z)$.  For any $T\subset [n]$, we denote by  $G^{(T)}$ the resolvent of $X^{(T)}$. When $T=\{i\}$ for some $i\in [n]$, we will write $G^{(i)}$ instead of $G^{(\{i\})}$, and for $T\subset [n]$ and $i \notin T$, $G^{(iT)}$ instead of $G^{(\{i\}\cup T)}$.

The Schur complement's formula \cite[(4.1)]{BGK} states that for any $i \in J^c$ and $z \in\mathbb{H}$,
\begin{equation} \label{Schur} \frac{1}{G_{ii}^{(J)}} =  -z - \sum_{k}^{(iJ)} X_{ik}^2 G_{kk}^{(iJ)} - \sum_{k\neq \ell}^{(iJ)} X_{ik} X_{i\ell} G_{k\ell}^{(iJ)},\end{equation}
where $\sum_k^{(T)} = \sum_{k\notin  T}$ for any $T\subset [n]$ and similarly for $\sum_{k\neq \ell}^{(T)}$. Very much in the spirit of the analytic proofs of local semicircle laws (see for example \cite{BGK}, \cite{EKHTY}), we will regard the system of equations \eqref{Schur} as a perturbation of the QVE defining $\upsilon_{W_n^{(J)}}$. Recall that $\upsilon_{W_n^{(J)}}$ is defined as the probability measure $\frac{1}{n} \sum_{i=1}^n \upsilon_i$, where $\upsilon_i$ has Stieltjes transform $m_i$ for any $i \in J^c$, and for any $z\in \mathbb{H}$,  $m(z)=(m_i(z))_{i\in J^c}$ is the unique solution  in $\mathbb{H}^{n}$ of the following QVE:
\begin{equation} \label{equpsilon} - \frac{1}{m_i(z)} = z + \sum_{j}^{(J)} X_{ij}^2 m_j(z), \ i\in J^c. \end{equation}
Showing that the equation \eqref{Schur} is a
small perturbation of the QVE \eqref{equpsilon} essentially amounts to prove that the cross terms in \eqref{Schur} vanish for most indices $i\in J^c$. While such a statement follows from a routine concentration argument in the proofs of local semicircle laws, the difficulty here is that we need such a statement to hold with overwhelming probability at the exponential scale $n^2p$. To this end, 
define for any $i\in J^c$ and $z\in \mathbb{H}$, the random variable
\begin{equation} \label{defYi} Y_{i}(z) := \sum_{k\neq \ell}^{(iJ)} X_{ik} X_{i\ell} G_{k\ell}^{(iJ)}.\end{equation} Note that although not explicitly specified in the notation, 
$Y_i(z)$ depends on $C$ which we regard here as a fixed parameter. With this notation, we will  prove the following proposition. 
\begin{Pro}\label{negli}
For any $\delta,\veps>0$, and $z\in \mathbb{H}$,
\[\lim_{n\to +\infty} \frac{1}{n^2p} \log \PP\big( \#\{ i\in J^c : |Y_{i}(z)|\geq \veps \} \geq \delta n\big)= - \infty.\]
\end{Pro}

This result is genuinely at the core of our strategy to understand the large deviation behaviour of the empirical spectral measure. We will need several intermediate steps to prove Proposition \ref{negli}. To give an idea why this result holds, we start by proving a tail bound for chaoses of order 2 of  ``sparse bounded random variables''. 

\begin{Lem}\label{chaos}Let $M$ be a $n\times n$ real symmetric matrix such that $\|M\|\leq 1$ and $(\sigma_k)_{k\in [n]}$ be a family of i.i.d. random variables with the same law as $\xi Z$, where $\xi$ is a Bernoulli random variable of parameter $p$, and $Z$ is a centered random variable bounded by $1$ independent of $\xi$. There exists $p_0\in(0,1)$ such that for any $p\leq p_0$ and $t>0$,
\[ \PP\Big( \big|\sum_{k\neq \ell}  M_{k\ell}\sigma_{k} \sigma_{\ell} \big|>t np \Big) \leq 2 e^{-np \widetilde{h}\big((t/16) \sqrt{\log(\frac{1}{p})}\big)},\]
where $\widetilde{h}(x) = \sup_{\theta \geq 0}\{ \theta x - \phi(\theta^2)\}$ for any $x\geq 0$, and $\phi(\theta) = e^\theta-1$ for any $\theta \geq0$. 
\end{Lem}
This simple tail bound explains why one can hope the proportion of non-trivial cross terms $Y_i(z)$ to be negligible at the exponential scale $n^2p$.  Indeed, conditionally on $J$ and $X^{(iJ)}$, $Y_i(z)$ can be written as a chaos of order $2$ in the  variables $\{ \xi_{ik} A_{ik}\}_{k \in J^c}$. As $\widetilde{h}(u) \sim_{+\infty} u \sqrt{\log u}$, it follows from Lemma \ref{chaos} that the tail distribution of $Y_i(z)$ decreases faster than $e^{-O(np \sqrt{\log (1/p) \log \log (1/p)})}$. If the variables $Y_i(z)$ were independent, this tail bound would immediately entails Proposition \ref{negli} using Bennett's inequality (see \cite[Theorem 2.9]{BLM}). Although the dependence of the variables $Y_i(z)$ poses some important difficulty, Lemma \ref{chaos} reveals how Proposition \ref{negli} can hold and will be a key element in its proof.
\begin{proof}
Let $\theta \geq 0$. We start by using a decoupling argument and introduce $(\widetilde{\sigma}_k)_{k\in [n]}$ an independent copy of $(\sigma_k)_{k\in [n]}$. By \cite[Theorem 3.1.1 (3.1.8)]{Pena}, we have that
\begin{equation} \label{logLaplace1} \EE\big( e^{\theta |\sum_{k\neq \ell} M_{k\ell} \sigma_k \sigma_\ell|} \big) \leq \EE\big( e^{8\theta | \sum_{k\neq \ell} M_{k\ell} \sigma_k \widetilde{\sigma}_\ell|}\big).\end{equation}
Let $u_\ell := \sum_{k\neq \ell} M_{k\ell} \sigma_k$ for any $\ell \in [n]$ and $u=(u_1,\ldots,u_n)^{\sf T}$. As $\|M\| \leq 1$, we have $\|u\|\leq 2\|\sigma\|$ where $\sigma = (\sigma_1,\ldots, \sigma_n)^{\sf T}$. Let $\Lambda$ be the log-Laplace transform of the joint law of the $\sigma_k$'s, that is $\Lambda(\zeta):= \log \EE (e^{\zeta \sigma_1})$, $\zeta\in \RR$. By \cite[(4.22)]{AB} we know that there exists $p_0\in (0,1)$ such that for $p\leq p_0$ and any $\zeta \in \RR$, $\Lambda(\zeta) \leq \zeta^2/\log(1/p)$.  
Assuming from now on that $p\leq p_0$ and denoting by $\widetilde{\sigma} = (\widetilde{\sigma}_1,\ldots,\widetilde{\sigma}_n)^{\sf T}$, it follows by independence that for any $\zeta \in \RR$,
\begin{equation} \label{logLaplace} \log \widetilde{\EE} \big( e^{\zeta \langle u,\widetilde{\sigma}\rangle}\big) \leq \frac{\zeta^2 \|u\|^2}{\log( \frac{1}{p})}\leq \frac{4\zeta^2 \|\sigma\|^2}{ \log( \frac{1}{p})},\end{equation}
where $\widetilde{\EE}$ denotes the expectation with respect to $\widetilde{\sigma}$. Now, using  \eqref{logLaplace}, we get 
\begin{equation} \label{logLaplace2} \EE\big( e^{8\theta | \sum_{k\neq \ell} M_{k\ell} \sigma_k \widetilde{\sigma}_\ell|}\big) \leq \EE\big(e^{8\theta \langle u,\widetilde{\sigma}\rangle}\big) + \EE\big(e^{-8\theta \langle u,\widetilde{\sigma}\rangle}\big) \leq 
2 \EE\Big[ e^{\frac{256 \theta^2 }{\log(\frac{1}{p})}\|\sigma\|^2} \Big].\end{equation}
Using the fact that $\sigma_k^2$ has the same distribution as $\xi Z^2$, $Z^2\leq 1$ and the concavity of the $\log$, we find that for any $\zeta \geq 0$ and $k\in [n]$,
\begin{equation} \label{logLaplace3} \log \EE(e^{\zeta \sigma_k^2}) \leq \log \EE(e^{\zeta \xi}) =\log (pe^\zeta +1-p) \leq p\phi(\zeta),\end{equation}
where $\phi(\zeta) = e^\zeta-1$. Using the independence of the $\sigma_k$'s and putting together \eqref{logLaplace1}, \eqref{logLaplace2} and \eqref{logLaplace3}, we have shown that for any $\theta \geq 0$, 
\[ \EE\big( e^{\theta |\sum_{k\neq \ell} M_{k\ell} \sigma_k \sigma_\ell|} \big) \leq 2e^{np \phi\big( \frac{256 \theta^2}{\log(1/p)}\big)}.\]
By Chernoff's inequality, this ends the proof of the claim.

\end{proof}

In order to circumvent the dependence between the variables $Y_i(z)$'s, we will use the following concentration inequality, which is a variation on a similar inequality due to Chatterjee \cite[Theorem 3.1]{Chtriangles}. 
\begin{Lem}\label{lemCha}
Let $(Z_i)_{i\in [n]}$ be a family of non negative random variables defined on the same probability space and $(\mathcal{F}^{(j)})_{j\in [n]}$ a family of sub-$\sigma$-algebras on that probability space. Denote by $\EE^{(j)}$ the conditional expectation given $\mathcal{F}^{(j)}$. Assume  that there exists $a>0$ and $\lambda_j>0$, $j\in [n]$ such that:
\begin{enumerate}
\item If $Z_j>0$, then $\sum_{i\in [n]} (Z_i - \EE^{(j)} (Z_i) )\leq a$, almost surely. 
\item For any $j\in [n]$, $ \EE^{(j)}(Z_j) \leq \lambda_j$, almost surely.
\end{enumerate}
Let $\lambda:=\sum_{j\in [n]} \lambda_j$. Then, for any $t> \lambda$, 
\[ \PP\big( \sum_{i\in [n]} Z_i>t\big) \leq \exp\Big(- \frac{\lambda}{a} h\Big( \frac{t}{\lambda}\Big)\Big)\leq \exp\Big( -\frac{t}{a}\log \frac{t}{3\lambda}\Big),\]
where $h(u) = u \log u - u+1$ for any $u\geq1$. 
\end{Lem}
This inequality can be seen as a generalisation of Bennett's inequality (see \cite[Theorem 2.9]{BLM} to dependent variables. Indeed, if the $Z_i$'s are independent Bernoulli variables and one takes $\mathcal{F}^{(j)}$ to be the $\sigma$-algebra generated by the variables $\{Z_j, j\neq i\}$, then $(1)$ holds with $a=1$ and $(2)$ with $\lambda_j = \EE (Z_j)$, so that one recovers Bennett's inequality. More generally, the relevant setup one should have in mind is when the variables $(Z_i)_{i\in[n]}$ are bounded functions of a family $(U_j)_{j\in [n]}$ of independent random variables in such a way that $Z_i$ depends ``mostly'' on $U_i$ for any $i \in [n]$. It is then natural to take  $\mathcal{F}^{(j)}$ as the $\sigma$-algebra generated by $\{U_i : i\neq j\}$, $j\in [n]$. In this situation, the parameter $a$ can be thought as a measure of the dependence between the $Z_i$'s: if  $Z_i$ depends ``mostly'' on $U_i$  then one expects that $Z_i \simeq \EE^{(j)} (Z_i)$ for any $i\neq j$, and one can hope  the sum in $(1)$ to be dominated by the $i=j$ term, in which case Lemma \ref{lemCha} yields that $\sum_{i\in [n]} Z_i$ behaves similarly as if the $Z_i$'s were independent.

\begin{proof}
Let $Z:= \sum_{i\in [n]} Z_i$ and denote by $\Lambda$ the the log-Laplace of $Z$ defined by $\Lambda(\theta):=\log \EE \exp(\theta Z)$ for any $\theta \in \RR$.  We will prove the following bound on the derivative of $\Lambda$. \begin{equation} \label{boundLambda} \Lambda'(\theta)\leq \lambda e^{\theta a}, \quad  \theta \geq 0.\end{equation}
Once this inequality proven, the claim will readily follow from \cite[Lemma 3.2]{Chtriangles} by integration and using Chernoff's inequality. Let $j\in [n]$ and $\theta \geq 0$. Using assumption (1) and the definition of the conditional expectation we get 
\begin{equation} \label{ineq1} \EE( Z_j e^{\theta Z}) \leq  \EE \big(Z_j e^{\theta \sum_{i\in [n]} \EE^{(j)} Z_i }\big) e^{\theta a}= \EE \big( [\EE^{(j)} Z_j]e^{\theta \sum_{i\in[n]} \EE^{(j)} Z_i} \big)e^{\theta a}.\end{equation}
Now, by assumption (2) and by conditional Jensen's inequality we obtain
\begin{equation} \label{ineq2} \EE \big( [\EE^{(j)} Z_j]e^{\theta \sum_{i\in[n]} \EE^{(j)} Z_i} \big) \leq \lambda_j \EE\big (e^{\theta \sum_{i\in[n]} \EE^{(j)} Z_i} \big)\leq \lambda_j \EE \big( e^{\theta Z}\big).\end{equation}
Combining the above two inequalities \eqref{ineq1} and \eqref{ineq2}, it follows that $\EE(Z_j e^{\theta Z}) \leq \lambda_j \EE(e^{\theta Z})$. Summing over $j\in [n]$, we get that $\EE(Z e^{\theta Z}) \leq \lambda \EE(e^{\theta Z})$. Since $\Lambda'(\theta) = \EE (Z e^{\theta Z}) /\EE (e^{\theta Z})$ for any $\theta \geq 0$, this gives the claimed bound \eqref{boundLambda}.
\end{proof}

In order to apply Lemma \ref{lemCha} to our observables of interest $(Y_i(z))_{i\in J^c}$, $z\in \mathbb{H}$, (after conditioning on the value of $J$), and in particular to check condition (1) of Lemma \ref{lemCha}, we introduce for any $i,j \in J^c$, $i\neq j$, a new variable $Y_{i(j)}(z)$ which will be a proxy for the conditional expectation $\EE_{\mathcal{F}^{(j)}} (Y_i(z))$, where $\mathcal{F}^{(j)}$ is the $\sigma$-algebra generated by $X^{(jJ)}$. More precisely, define for any $i,j \notin {J}$, $i\neq j$ and $z\in \mathbb{H}$, 
\begin{equation} \label{defYij} Y_{i(j)}(z) := \sum_{k\neq \ell}^{(jJ)} X_{ik} X_{i\ell} G^{(ijJ)}_{k\ell}(z).\end{equation}
To ease the notation, we denote in the sequel by $\lesssim$ if an inequality holds up to a absolute multiplicative constant and $\lesssim_z$ if the inequality holds up to a multiplicative constant which depends on $z$.
 As preparatory work towards the proof of Proposition \ref{negli}, we derive the following bound on the difference between ${Y}_i$ and ${Y}_{i(j)}$, which we will use to check assumption (1) of Lemma \ref{lemCha}.
 \begin{Lem}
\label{diff}
For any $i,j \in J^c$, $i\neq j$, and $z\in \mathbb{H}$,
\[ |{Y}_i(z)-{Y}_{i(j)}(z)|\lesssim_z (1+\|X_j\|^4+\|X_i\|^{12})\big[X_{ij}^2+ |G_{ij}^{(J)}|^2 + \sum_{k}^{(J)} X_{ik}^2|G_{kj}^{(J)}|^2\big],\]
where $X_\ell$ denotes the $\ell^{\text{th}}$-column of $X$ for any $\ell \in [n]$.
\end{Lem}

\begin{proof}
We will make use repetitively of the following two resolvent identities (see \cite[(3.4) and (3.5)]{BGK}). The first one compares the entries of $G^{(T)}$ and $G^{(jT)}$ for some $T\subset [n]$ and $j \notin T$:
\begin{equation} 
\label{resolvid1}
G_{k \ell}^{(T)} = G_{k \ell}^{(jT)} + \frac{G_{kj}^{(T)} G_{j\ell}^{(T)} }{G_{jj}^{(T)}}, \quad  k,\ell \notin T\cup \{j\}.
\end{equation}
The second resolvent identity gives a sort of recursive formula for the off-diagonal entries of the resolvent $G_{ij}^{(T)}$ for some $T\subset [n]$, and $i,j \notin T$, $i\neq j$,
\begin{equation} 
\label{resolvid2}
G_{i j}^{(T)} = - G_{ii}^{(T)} \sum_{k}^{(iT)} X_{i k} G_{k j}^{(iT)} = -G_{jj}^{(T)} \sum_{\ell}^{(jT)}G_{i \ell}^{(jT)} X_{\ell j}.
\end{equation}
Further, we note that for any $T\subset [n]$ and $k,\ell \notin T$, we have $G^{(T)}_{k\ell} = G^{(T)}_{\ell k}$ as $X$ is symmetric. 
Now, fix $i,j \in J^c$ such that $i\neq j$. Splitting the sum defining $Y_i(z)$ \eqref{defYi} into two parts, the first part containing the indices $k,\ell \neq j$ and the second where either $k$ or $\ell$ is equal to $j$, we define 
\[ \Delta_1 :=   \sum_{k\neq \ell}^{(ijJ)} X_{ik}X_{i\ell} (G_{k\ell}^{(iJ)} - G_{k\ell}^{(ijJ)}), \ \Delta_2 = 2 X_{ij} \sum_{\ell}^{(ijJ)} X_{i\ell} G_{j\ell}^{(iJ)}.\]
With this notation, we have $Y_i(z) - Y_{i(j)}(z)=\Delta_1+\Delta_2$.
We will bound separately $\Delta_1$ and $\Delta_2$ using the resolvent identities \eqref{resolvid1} and \eqref{resolvid2}. 
Starting with $\Delta_1$, we get using \eqref{resolvid1} that
\begin{align*}
\Delta_1 &= \frac{1}{G_{jj}^{(iJ)}}\sum_{k\neq \ell}^{(ijJ)} X_{ik}X_{i\ell} G_{kj}^{(iJ)} G_{j\ell}^{(iJ)} \\
&= \frac{1}{G_{jj}^{(iJ)}} \Big[\Big( \sum_k^{(iJ)} X_{ik} G_{kj}^{(iJ)}\Big)^2-\sum_{k}^{(iJ)} X_{ik}^2(G_{kj}^{(iJ)})^2-2X_{ij} G_{jj}^{(iJ)}\sum_k^{(ijJ)} X_{ik} G_{kj}^{(iJ)}\Big].
\end{align*}
Together with the identity \eqref{resolvid2}, this implies that
\begin{equation} \label{decom}  \Delta_1 = \frac{1}{G_{jj}^{(iJ)}}\Big[ \Big( \frac{G_{ij}^{(J)}}{G_{ii}^{(J)}}\Big)^2 -\sum_{k}^{(iJ)} X_{ik}^2(G_{kj}^{(iJ)})^2  +2X_{ij} G_{jj}^{(iJ)}\frac{G_{ij}^{(J)}}{G_{ii}^{(J)}} + 2X_{ij}^2(G_{jj}^{(iJ)})^2  \Big].\end{equation}
Note that by Schur's complement formula  \cite[(4.1)]{BGK}, for any $T\subset [n]$ and $i\notin T$, 
\begin{equation} \label{boundresolv}\frac{1}{|G_{ii}^{(T)}|} = |z + \langle X_i^{(T)}, G^{(T)} X_i^{(T)}\rangle |\leq |z|+\frac{\|X_i^{(T)}\|^2}{\Im z},\end{equation}
where we used the fact that $\|G^{(T)}\|\leq 1/\Im z$.  Thus, $1/|G_{ii}^{(T)}|\lesssim_z \|X_i\|^2 + 1$ for any $i\notin T$. Making use of this bound in \eqref{decom} and the fact that $|G_{jj}^{(iJ)}|\lesssim_z 1$, we get
\begin{align}
|\Delta_1|&\lesssim_z  (\|X_j\|^4+\|X_i\|^8+1)|G_{ij}^{(J)}|^2 +(\|X_j\|^2+1)\sum_k^{(iJ)} X_{ik}^2 |G_{kj}^{(iJ)}|^2 \nonumber \\
&+ (\|X_i\|^2+1)|X_{ij}G_{ij}^{(J)}| + X_{ij}^2,\label{decom1}
\end{align}
where we used the inequality $(\|X_j\|^2+1)(\|X_i\|^2+1)^2\lesssim \|X_j\|^4+\|X_i\|^8+1$.
By the resolvent identity \eqref{resolvid1}, we find that 
\[ \sum_k^{(iJ)} X_{ik}^2 |G_{kj}^{(iJ)}|^2\lesssim \sum_k^{(iJ)} X_{ik}^2 |G_{kj}^{(J)}|^2 + \sum_k^{(J)} X_{ik}^2 \frac{|G_{ik}^{(J)} G_{ij}^{(J)}|^2}{|G_{ii}^{(J)}|^2}.\]
Using again \eqref{boundresolv} and the bound $|G_{ik}^{(J)}|\lesssim_z 1$ for any $k \notin J$, this yields
\[ \sum_k^{(iJ)} X_{ik}^2 |G_{kj}^{(iJ)}|^2\lesssim_z \sum_k^{(iJ)} X_{ik}^2 |G_{kj}^{(J)}|^2 +(\|X_i\|^6+1) |G_{ij}^{(J)}|^2,\]
where we used that $\|X_i\|^2(\|X_i\|^2+1)^2\lesssim \|X_i\|^6+1$.
Plugging this estimate in \eqref{decom1} and using the inequality $|X_{ij}G_{ij}^{(J)}|\leq 2(|X_{ij}|^2+|G_{ij}^{(J)}|^2)$, it follows that 
\[ |\Delta_1| \lesssim_z (1+\|X_j\|^4+\|X_i\|^{12}) \big[ |G_{ij}^{(J)}|^2 + \sum_{k}^{(J)} X_{ik}^2|G_{kj}^{(J)}|^2+X_{ij}^2\big]\]
Finally, we bound $\Delta_2$. Using \eqref{resolvid1} we find that $\Delta_2 = -2X_{ij} G_{ij}^{(J)}/G_{ii}^{(J)} - 2X_{ij}^2 G_{jj}^{(iJ)}$. Applying the bound \eqref{boundresolv}, the inequality $|X_{ij}G_{ij}^{(J)}|\lesssim X_{ij}^2+ |G_{ij}^{(J)}|^2$ and $|G_{jj}^{(iJ)}|\lesssim_z 1$, we get 
\[ |\Delta_2|\lesssim_z (1+\|X_i\|^2) (X_{ij}^2+|G_{ij}^{(J)}|^2).\]
As $Y_i(z)-Y_{i(j)}(z) = \Delta_1+\Delta_2$, this ends the proof of the claim. 
\end{proof}

We are now ready to prove Proposition \ref{negli}.
\begin{proof}[Proof of Proposition \ref{negli}]
Let $\veps>0$, $z\in \mathbb{H}$ and $C\geq 1$ be large enough.  Let $\phi_\veps$ be a $2/\veps$-Lipschitz function such that ${\bf 1}_{[\veps,+\infty)}\leq \phi_\veps \leq {\bf 1}_{[\veps/2,+\infty)}$ and define for any $i \in J^c$,  $Z_i := \phi_\veps(|Y_{i}(z)|\big)$, where $Y_i(z)$ is defined in \eqref{defYi}.
For any $\mathcal{J} \subset [n]$ we denote by $\PP_{\mathcal{J}}$ the conditional probability measure $\PP( . \mid J= \mathcal{J})$. 
We will prove that for any $\delta>0$, \begin{equation} \label{claim} \lim_{n\to+\infty} \frac{1}{n^2p}\sup_{\mathcal{J}\subset [n]} \log \PP_\mathcal{J} \big( \sum_{i\notin \mathcal{J}} Z_i > \delta n\big) = -\infty.\end{equation}
Assume for the moment that \eqref{claim} holds. Note that $\PP( \#\{ i\in J^c : |Y_i(z)|\geq \veps\} \geq \delta n ) \leq \PP( \sum_{i\notin J} Z_i \geq \delta n)$. As the number of subsets of $[n]$ is $2^n$ and $np \gg 1$, a union bound together with \eqref{claim}  yields indeed the claim. 
 We now move on to prove \eqref{claim}. To this end, fix $\delta>0$ and $\mathcal{J} \subset [n]$. Aiming at applying Lemma \ref{lemCha}, we naturally consider for any $j\notin \mathcal{J}$, $\mathcal{F}^{(j)}$ to be the $\sigma$-algebra generated by $X^{(j\mathcal{J})}$. In order to check condition (1) of Lemma \ref{lemCha}, we introduce or any $i\notin J$, $j\neq i$,  the variable $Z_{i(j)} :=\phi_\veps(Y_{i(j)}(z))$, where $Y_{i(j)}(z)$ is defined in \eqref{defYij}.
We claim that almost surely under $\PP_{\mathcal{J}}$, for any $j\notin \mathcal{J}$, 
\begin{equation} \label{claimbound} \sum_{i}^{(\mathcal{J})} |Z_i-Z_{i(j)}|\lesssim_{C,z,\veps} 1,\end{equation}
where $\lesssim_{C,z,\veps}$ denotes that the inequality holds up to a multiplicative constant depending on $C,z,\veps$. 
Fix $j \in \mathcal{J}^c$. Using the fact that $\phi_\veps$ is $2/\veps$-Lipschitz and Lemma \ref{diff}, we deduce that for any $i\in \mathcal{J}^c$, $i\neq j$, on the event $\{J=\mathcal{J}\}$,
\[ |Z_i-Z_{i(j)}| \lesssim_z 2\veps^{-1} (1+\|X_j\|^4+\|X_i\|^{12}) \big[X_{ij}^2+|G_{ij}^{(\mathcal{J})}|^2 + \sum_{k}^{(\mathcal{J})} X_{ik}^2|G_{jk}^{(\mathcal{J})}|^2\big].\]
Summing over $i\in \mathcal{J}^c$, $i\neq j$, and using the fact that for any $k\in {J}^c$, $\|X_k\|\leq C$ by definition of $J$ (see \eqref{defJ}), it follows that $\PP_\mathcal{J}$-almost surely, 
 \begin{align*} \sum_{i}^{(j\mathcal{J})}|Z_i-Z_{i(j)}|&\lesssim_{C,z,\veps}   \sum_{i}^{(\mathcal{J})} |G_{ij}^{(\mathcal{J})}|^2+\|X_j\|^2+\sum_{k}^{(\mathcal{J})}\|X_k\|^2 |G_{jk}^{(\mathcal{J})}|^2 \\
 &\lesssim_{C,z,\veps} \sum_{k}^{(\mathcal{J})}|G_{jk}^{(\mathcal{J})}|^2 +1. \end{align*}
By Ward identity \cite[(3.6)]{BGK},  we have $\sum_{k \notin \mathcal{J}} |G_{jk}^{(\mathcal{J})}|^2 = \Im G_{jj}^{(\mathcal{J})}/\Im z \leq 1/(\Im z)^2$, which ends the proof of \eqref{claimbound}. 

To ease the notation, denote by $\EE_\mathcal{J}^{(j)}$  the conditional expectation under $\PP_{\mathcal{J}}$ given $\mathcal{F}^{(j)}$. Denote for any $j \notin \mathcal{J}$ by $\EE^{(j)}$ the conditional expectation given $\mathcal{F}^{(j)}$ under $\PP$. Note that almost surely under $\PP_\mathcal{J}$, $Z_{i(j)}$  is a measurable function of $X^{(j\mathcal{J})}$ for $i,j \in \mathcal{J}^c$, $i\neq j$.  Therefore $\EE_{\mathcal{J}}^{(j)} (Z_{i(j)}) = Z_{i(j)}$ for any $i,j \in \mathcal{J}^c$, $i\neq j$. It follows by taking the conditional expectation given $\mathcal{F}^{(j)}$ under $\PP_{\mathcal{J}}$ in \eqref{claimbound} that $\sum_{i}^{(j\mathcal{J})} |\EE_\mathcal{J}^{(j)} Z_i-Z_{i(j)} |\lesssim_{C,\kappa,\veps} 1$, almost surely under $\PP_{\mathcal{J}}$. Together with \eqref{claimbound}, this entails that 
\begin{equation} \label{cond1} \sum_{i}^{(j\mathcal{J})} |Z_i- \EE^{(j)} Z_i |\lesssim_{C,z,\veps} 1, \quad \text{ a.s. under $\PP_\mathcal{J}$}.\end{equation}
Now, since $|Z_j|\leq 1$, the above bound \eqref{cond1} is unchanged even by adding the $i=j$ term, so that the first condition of Lemma \ref{lemCha} holds with $a$ equal to a certain positive constant depending on the parameters $C,z$ and $\veps$. 

Moving on to check the second assumption of Lemma \ref{lemCha}, we fix $j\in \mathcal{J}^c$ and compute $\EE^{(j)}_{\mathcal{J}}(Z_j)$. Denote by $\EE^{(j)}$ the conditional expectation given $\mathcal{F}^{(j)}$ under $\PP$. Observe that $\EE^{(j)}_{\mathcal{J}}(Z_j)  \leq \PP^{(j)}_{\mathcal{J}}(|Y_j|\geq \veps/2)$ almost surely since $\phi_\veps \leq \Car_{[\veps/2,+\infty)}$. We claim that up to paying a factor $e^{-O(np)}$ we can remove the conditioning given $J=\mathcal{J}$, meaning that 
\begin{equation}\label{boundcond}  \PP^{(j)}_{\mathcal{J}}(|Y_j(z)|\geq \veps/2) \leq (1-p)^{-n} \PP^{(j)}(|Y_j(z)|\geq \veps/2), \ \PP_\mathcal{J}\text{-a.s.}\end{equation}
We will show more generally that for any non-negative $X^{(\mathcal{J})}$-measurable random variable $T$, 
\begin{equation} \label{claimgene}\EE_{\mathcal{J}}^{(j)} (T)\leq (1-p)^{-n}\EE^{(j)}(T), \ \PP_\mathcal{J}\text{-a.s.}\end{equation}
This will prove the claim \eqref{boundcond} as $Z_j$ is $\PP_\mathcal{J}$-almost surely equal to a  $X^{(\mathcal{J})}$-measurable random variable.
Observe first that denoting by $E$ the event $\{ \forall i\in \mathcal{J}^c, i\neq j, \|X_i^{(j)}\|\leq C, \forall i\in \mathcal{J}, \|X_i\|>C\}$, we have the inequalities
\begin{equation} \label{encard} {\bf 1}_{E\cap \{\forall i \in \mathcal{J}^c, \xi_{ij}=0\}}\leq  {\bf 1}_{\{ J=\mathcal{J} \}} \leq {\bf 1}_{E}.\end{equation}
Moreover, $E$ is measurable with respect to the variables $(X_i)_{i\in \mathcal{J}}$ and $X^{(j\mathcal{J})}$. It follows that if $T$ is $X^{(\mathcal{J})}$-measurable, then $\Car_{E}$ and $T$ are conditionally independent given $X^{(j\mathcal{J})}$. 
Denote by $\mathcal{H}_{(j\mathcal{J})}$ the set of $(\mathcal{J}\cup \{j\})^c \times (\mathcal{J}\cup\{j\})^c$ symmetric matrices. Thus, for any non negative measurable function $f:\mathcal{H}_{(j\mathcal{J})} \to \RR_+$ we can write using the fact that $T$ is also non-negative and \eqref{encard},
\begin{align}
\EE\big( Tf(X^{(j\mathcal{J})}) {\bf 1}_{\{J= \mathcal{J}\}}\big)&\leq  \EE\big[ T f(X^{(j\mathcal{J})}) \Car_{E}\big] =\EE\Big[\EE^{(j)} (T) f(X^{(j\mathcal{J})}) \PP^{(j)}(E) \Big] \nonumber\\
& = \EE\Big[\EE^{(j)} (T) f(X^{(j\mathcal{J})}) \Car_E \Big]. \label{maj}
\end{align}
Now, since $\{\forall i\in \mathcal{J}^c, \xi_{ij}=0\}$ is independent from $((X_i)_{i\in \mathcal{J}}, X^{(j\mathcal{J})})$, we deduce using \eqref{encard} that 
\begin{align} \EE\big[\EE^{(j)} (T) f(X^{(j\mathcal{J})}) \Car_{J=\mathcal{J}} \big]&\geq \EE\big[\EE^{(j)} (T) f(X^{(j\mathcal{J})})  {\bf 1}_{E\cap \{\forall i \in \mathcal{J}^c, \xi_{ij}=0\}}  \big]  \nonumber \\
& \geq (1-p)^{|\mathcal{J}^c|}  \EE\big[\EE^{(j)} (T) f(X^{(j\mathcal{J})}) \Car_E \big].\label{min}
\end{align}
Putting together \eqref{maj} and \eqref{min}, we have shown that for any non-negative measurable function $f : \mathcal{H}_{(j\mathcal{J})} \to\RR_+$, 
\[\EE_\mathcal{J}[Tf(X^{(j\mathcal{J})})] \leq (1-p)^{-|\mathcal{J}^c|}\EE_\mathcal{J}[\EE^{(j)} (T)f(X^{(j\mathcal{J})})],\] which implies \eqref{claimgene}.

Finally, using that $\|G^{(j\mathcal{J})}\|\leq 1/\Im z$ and the fact that for any $k,\ell\in [n]$, $|A_{k\ell}|\leq R$ a.s., it follows from Lemma \ref{chaos} that  
\[ \PP^{(j)} (|Y_j(z)|\geq \veps/2) \leq e^{-np \alpha_n}, \ \text{a.s.}\]
where $\alpha_n:=\alpha_n(R,\veps,z)\to +\infty$ as $n\to +\infty$. Using \eqref{boundcond}, we deduce that $\PP_\mathcal{J}$-almost surely,  $\EE^{(j)}_{\mathcal{J}}(Z_j) \leq (1-p)^n e^{-np \alpha_n}$. Thus, the second condition of Lemma \ref{lemCha} holds with  $\lambda =  n(1-p)^n e^{-np \alpha_n}$. As $\alpha_n\gg1$ and $np \gg \log n$, Lemma \ref{lemCha} implies that for $n$ large enough
\[ \PP_\mathcal{J}\big(\sum_{i\notin \mathcal{J}} Z_i >\delta n  \big) \leq \exp\Big(-\frac{\delta n}{a}\log \frac{\delta n}{3\lambda}\Big)\leq  \exp\Big(-\frac{\delta}{2a} n^2 p \alpha_n\Big).\]
As $a$ is some constant depending on $C,z,\veps$, and $\alpha_n\to +\infty$ as $n\to+\infty$, this ends the proof.

\end{proof}

\subsection{Stability of QVEs}
To leverage the result of Proposition \ref{negli}, we develop in this section a stability estimate for solutions of Quadratic Vector Equations that is suited for our purpose, based on the results of Ajanki, Erd\H{o}s and Kruger \cite{AEK}.

To this end, we are considering a more general form of QVEs. Following the notation of \cite{AEK}, consider a set of label $\mathfrak{X}$ and $\mathcal{B}$ the set of bounded complex-valued functions on $\mathfrak{X}$, that is, $\mathcal{B}:=\{ w \in \CC^{\mathfrak{X}} : \sup_{x\in \mathfrak{X}} |w_x|<+\infty\}$. Typically, $\mathfrak{X}$ will be taken in the applications either as a finite set $\{1,\ldots,n\}$ or as $[0,1]$. We endow $\mathcal{B}$ with the sup norm $\|\ \|_\infty$, defined by $\|w\|_\infty= \sup_{x\in \mathfrak{X}} |w_x|$ for any $w\in \CC^{\mathfrak{X}}$, thus making $\mathcal{B}$ a Banach space. We moreover consider a probability measure $\pi$ on  $\mathfrak{X}$ (which will be either the uniform measure if $\mathfrak{X}$ is finite or the Lebesgue measure if $\mathfrak{X}=[0,1])$, and denote by $\langle . , .\rangle$ the inner product on the complex vector space $L^2(\mathfrak{X},\pi)$, and by $\| \ \|_{L^2}$ the associated $L^2$ norm. 

Now, for an operator $S :\mathcal{B} \to \mathcal{B}$, we denote by $\|S\|_\infty$ its operator norm with respect to the sup norm on $\mathcal{B}$. We say that $S$ is {\em bounded} if $\|S\|_\infty<+\infty$, and that $S$  is {\em symmetric and positivity preserving} if  for any $u,v\in \mathcal{B}$, and non negative $p\in \mathcal{B}$,
\[ \langle u, Sv\rangle = \langle v,Su\rangle, \quad  (Sp)_x\geq 0, \ \forall x\in \mathfrak{X}.\]
We denote moreover by $\mathcal{B}^+$ the subset of $\mathcal{B}$ consisting of functions $w$ such that $\Im w > 0$. With this notation, we can now state the following stability estimate.

\begin{Lem}[Stability]\label{stability}Let $S :\mathcal{B} \to \mathcal{B}$ be a bounded symmetric and positivity preserving operator, and let for any $z\in \mathbb{H}$, $d(z)\in \mathcal{B}$.
Assume that for any $z\in \mathbb{H}$, $\widetilde{m}(z)=(\widetilde{m}_x(z))_{x\in\mathfrak{X}}$ satisfies  the equation
\[ - \frac{1}{\widetilde{m}(z)} = z + S\widetilde{m}(z) + d(z),\]
For any $z\in\mathbb{H}$, let  ${m}(z) = ({m}_x(z))_{x\in\mathfrak{X}}$ be the unique solution in $\mathcal{B}^+$ of the QVE
\begin{equation} \label{QVES} - \frac{1}{{m}(z)} = z+ S{m}(z).\end{equation}
Then, there exists a numerical constant $\kappa \geq 1$ such that if $\Im z \geq [ \kappa (\|S\|_\infty \vee 1)^2]\vee |\Re z|$, then
\[  \| m(z) - \widetilde{m}(z)\|_{L_2} \leq \kappa \big( \|S\|_\infty\vee 1) \|d\|_{L_2}.\]
\end{Lem}

Our stability estimate differs from the one of \cite[Theorem 2.12]{AEK} in terms of both assumptions on $S$ and control on the distance between the solutions. While Lemma \ref{stability} only assumes that the operator $S$ is symmetric and positivity preserving, \cite[Theorem 2.12]{AEK} works under stronger smoothing and primitivity assumptions (see \cite[\textbf{A2-3}]{AEK}), which we cannot afford in our setting. On the other hand, our estimate is much coarser, as it only holds in $L^2$ norm instead of sup norm and the solutions are evaluated far from the real axis, which greatly simplifies the arguments.

Technically, the argument to prove Lemma \ref{stability} is only a simpler variation of the one of \cite[Theorem 2.12]{AEK}. As shown in \cite{AEK}, the stability of the QVE \eqref{QVES} is related to the inverse of a certain auxiliary operator $B$ and the estimation of the norm of this inverse. More precisely, for $z\in \mathbb{H}$, let $F : \mathcal{B} \to \mathcal{B}$ be the operator defined by 
\[ Fw = |m(z)| S(|m(z)|w), \quad w \in \mathcal{B},\]
  where $m(z)$ is the unique solution of \eqref{QVES} in $\mathcal{B}^+$. Note that although not explicitly mentioned in the notation, $F$ depends on $z$. Since $\|m(z)\|_\infty \leq (\Im z)^{-1}$ by \cite[Theorem 2.1]{AEK} and $S$ is symmetric and bounded, $F$ is clearly also symmetric and bounded. Now, define $B : \mathcal{B} \to \mathcal{B}$ as $B = e^{-2iq(z)} -F$, where $q(z) : \mathfrak{X} \to [0,2\pi)$ is a representation of the argument of $m(z)$, that is such that $e^{iq(z)} = m(z)/|m(z)|$. Before going into the proof of Lemma \ref{stability}, we repeat some preliminary discussion from \cite[section §4.2]{AEK}. Note that since $S$ is symmetric, $Sw=0$ $\pi$-a.s. whenever $w=0$ $\pi$-a.s., so that $S$ is well-defined as a bounded operator on $L^\infty(\mathfrak{X},\pi)$ and $\|S\|_{L^\infty\to L^\infty}\leq \|S\|_\infty$. Using the duality of the $L^p$ norms and the symmetry of $S$, one can show that $S$ extends as a bounded operator on $L^1(\mathfrak{X},\pi)$ and that $\|S\|_{L^1\to L^1} = \|S\|_{L^\infty \to L^\infty}$. The Riesz-Thorin interpolation theorem entails that $S$ is a bounded operator on $L^2(\mathfrak{X},\pi)$, and that $\|S\|_{L^2\to L^2} \leq \sqrt{\|S\|_{L^1\to L^1} \|F\|_{L^\infty\to L^\infty}}\leq\|S\|_\infty$. Since $F$ is also a bounded symmetric operator, the same holds for $F$ as well.
 In our setting, the solutions of the QVEs are evaluated far from the real axis, so that the norm of $B^{-1}$, as an operator on $L^2(\mathfrak{X},\pi)$, can be easily bounded as follow.

 \begin{Lem}\label{boundnormB}
 For any $z \in \mathbb{C}$, $\Im z\geq 1$, $\|B^{-1}\|_{L^2 \to L^2} \leq 16(\|S\|_\infty\vee1) |z|^2/(\Im z)^2$.  
 \end{Lem}
 
 \begin{proof}
 To prove the invertibility of $B$ as an operator on $L^2(\mathfrak{X},\pi)$ and the claimed bound on its inverse, we will prove that $\|F\|_{L^2\to L^2} \leq 1-\eta$, where $\eta := \frac{(\Im z)^2}{16(\|S\|_\infty\vee 1)|z|^2}$. Once this bound proven, it follows readily that $B$ as an operator on $L^2(\mathfrak{X},\pi)$ is invertible with $B^{-1} = \sum_{k\in \NN} F^k e^{-(k+1)2iq}$, and as a consequence $\|B^{-1}\|_{L_2\to L_2}\leq \eta^{-1}$. It was shown in the proof of \cite[Lemma 4.5 (4.36)]{AEK} that 
 \begin{equation} \label{boundL2F} \|F\|_{L^2\to L^2} \leq 1- \veps, \quad \veps:= \Im z \inf_{x\in\mathfrak{X}} \frac{|m_x(z)|^2}{\Im m_x(z)} \in (0,1].\end{equation}
This bound was proven by exhibiting in  \cite[(4.36)]{AEK}, a vector $w\in L^2(\mathfrak{X},\pi)$ positive almost surely such that $Fw \leq (1-\veps) w$, and by using a Perron-Frobenius type theorem (see \cite[Lemma 4.6]{AEK}) to infer the bound \eqref{boundL2F} on the operator norm.
Now, by \cite[Theorem 2.1]{AEK} we know that for any $x\in \mathfrak{X}$, $m_x$ can be written as the Stieltjes transform of a probability measure $v_x$ on $\RR$ which is  supported on the interval $[-2\|S\|_\infty^{1/2},2\|S\|_\infty^{1/2}]$. Thus, it follows that
 \[ \Im m_x(z) = \int \frac{\Im z }{(u-\Re z)^2+(\Im z)^2} d\upsilon_x(u) \geq \frac{ \Im z}{(2\|S\|_\infty^{1/2} +|z|)^2}\geq \frac{\Im z }{8(\|S\|_\infty+|z|^2)}.\]
 Now, using the bound $|m_x(z)|^2\geq (\Im m_x(z))^2$, we obtain that $\veps \geq \eta$, which ends the proof. 
  \end{proof}

We are now ready to prove Lemma \ref{stability}.
\begin{proof}[Proof of Lemma \ref{stability}]Let $z\in \CC$ such that $\Im z \geq [128 (\|S\|_\infty\vee1)^2] \vee |\Re z|$. For sake of clarity, we drop the $z$-dependence in our notation. 
We repeat the proof of \cite[Lemma 5.11]{AEK} to relate the difference between $m(z)$ and $\widetilde{m}(z)$ to the $L^2$ operator norm of $B^{-1}$. Let $h = m -\widetilde{m}$. By \cite[Lemma 5.11 (5.60)]{AEK}, we know that $h$ satisfies the equation
\[ h = |m| B^{-1}[ e^{-iq} hSh +(|m|+e^{-iq}h)d].\]
Taking $L^2$-norms, using the fact that $\|S\|_{L_2\to L_2}\leq \|S\|_\infty$ and the bound $\|m\|_\infty\leq 1$ as $\Im z \geq 1$,  we get
\begin{align}
 \|h\|_{L^2}&\leq \|m\|_\infty \|B^{-1}\|_{L^2\to L^2}  \|h\|_\infty  \|S\|_\infty \|h\|_{L^2}+ \|m\|_\infty^2 \|B^{-1}\|_{L^2\to L^2} \|d\|_{L^2}\nonumber \\
&+ \|m\|_\infty\|B^{-1}\|_{L^2\to L^2} \|h\|_\infty \|d\|_{L^2}\nonumber \\
& \leq \big[ \|B^{-1}\|_{L^2\to L^2}\|S\|_\infty\|h\|_\infty] \|h\|_{L^2} + \|B^{-1}\|_{L^2\to L^2} [1+\|h\|_\infty]\|d\|_{L^2}.\label{ineqnorm}
\end{align}
Now, using the fact that $\|h\|_\infty \leq \|m\|_\infty+\|\widetilde{m}\|_\infty\leq  2/\Im z$ and Lemma \ref{boundnormB}, we find that
\[ \|B^{-1}\|_{L^2\to L^2} \|S\|_\infty \|h\|_\infty\leq 32 (\|S\|_\infty\vee 1)^2 \frac{|z|^2}{(\Im z)^3}\leq  64 \frac{(\|S\|_\infty\vee 1)^2}{\Im z} \leq 1/2,\]
where we used the fact $\Im z \geq [128( \|S\|_\infty \vee 1)^2]\vee  |\Re z|$. Using  again the bound $\|h\|_\infty\leq 2/\Im z \leq 1$, we find by coming back to \eqref{ineqnorm}, that 
\[ \|h\|_{L^2} \leq 4\|B^{-1}\|_{L^2\to L^2}\|d\|_{L^2} \leq 64 (\|S\|_\infty\vee 1) \frac{|z|^2}{(\Im z)^2} \|d\|_{L^2} \leq 128 (\|S\|_\infty \vee 1)\|d\|_{L^2},\]
where we used again the fact that $|\Re z |\leq \Im z$. 
\end{proof}

Using the stability estimate of Lemma \ref{stability} and Proposition \ref{negli} we can finally prove the following exponential approximation. 
\begin{Lem}\label{equivexpo}
For any $\veps>0$,
\[ \lim_{C\to +\infty}\limsup_{n\to+\infty} \frac{1}{n^2p} \log \PP\big( d(\mu_{X^{(J_C)}},\upsilon_{W_n^{(J_C)}})>\veps\big) =-\infty,\]
where $d$ is defined in \eqref{defd}.
\end{Lem}

Once this lemma proven, the conclusion of Proposition \ref{equivexpomain1} immediately follows by putting together Lemmas \ref{equivexpo}, \ref{approxW} and \ref{reduc}. 
\begin{proof}We start by a series of reductions of the statement to the case of estimating the probability that the Stieltjes transforms of $\mu_{X^{(J)}}$ and $\upsilon_{W_n^{(J)}}$ are pointwise apart. 
Using Lemma \ref{deg}, we see that it suffices to prove that for any $\delta>0$ and $C\geq 1$, 
\begin{equation} \label{reduc2} \lim_{n\to+\infty} \frac{1}{n^2p} \log \PP\big( d(\mu_{X^{(J)}},\upsilon_{W_n^{(J)}})>\delta, \ |J|\leq n/2\big) =-\infty.\end{equation}
Define for any $\mu,\nu \in \mathcal{P}(\RR)$, 
\[ d_C(\mu,\nu) := \sup\big\{|m_\mu(z)-m_\nu(z)| :  z\in \mathcal{D}_C\big\},\] 
where $\mathcal{D}_C:=\{z\in \mathbb{H} : |\Re z|\leq \kappa C^4 \leq \Im z \leq 2\kappa C^4\}$, and $\kappa \geq 1$ is the numerical constant in Lemma \ref{stability}. Since $\mathcal{D}_C$ contains an accumulation point in $\mathbb{H}$, $d_C$ defines a distance on $\mathcal{M}_1(\RR)$, the set of Borel measures on $\RR$ with total mass less than $1$, which is compatible with the vague topology. Clearly, the same holds for $d$. Thus, $d$ and $d_C$ are two equivalent distances on $\mathcal{M}_1(\RR)$ endowed with the vague topology, which is a compact topological space by Helly's Selection Theorem. As a result, $d$ and $d_C$ are uniformly equivalent metrics, that is, for any $\delta>0$, there exists $\delta'>0$ such that for any $\mu,\nu \in \mathcal{M}_1(\RR)$, if $d_C(\mu,\nu) \leq \delta'$ then $d(\mu,\nu) \leq \delta$. We deduce that in order to prove \eqref{reduc2}, it is enough to show that for any $C\geq 1$ and $\delta'>0$, 
\begin{equation} \label{reduc3} \lim_{n\to+\infty} \frac{1}{n^2p} \log \PP\big( d_C(\mu_{X^{(J)}},\upsilon_{W_n^{(J)}})>\delta', |J|\leq n/2\big) =-\infty.\end{equation}
Denote by $\widetilde{g}$ the Stieltjes transform of $\mu_{X^{(J)}}$ and by $g$ the one of $\upsilon_{W_n^{(J)}}$.
Note that $z \mapsto m_\mu(z)$ is $1$-Lipschitz on $\{z \in \mathbb{H} : \Im z \geq 1\}$ for any probability measure $\mu\in\mathcal{P}(\RR)$. Using a union bound we see that it is sufficient to show that for any $C\geq 1$, $z\in \mathcal{D}_C$, and any $\delta'>0$,
\begin{equation} \label{reduc4} \lim_{n\to+\infty} \frac{1}{n^2p} \log \PP\big( |\widetilde{g}(z) - g(z)|>\delta', |J|\leq  n/2\big) =-\infty.\end{equation}
Fix $C\geq 1$ and $z\in \mathcal{D}_C$.
Now, set for any $i\in J^c$, $\widetilde{m}_i(z) := G_{ii}^{(J)}(z)$. By the Schur complement formula \eqref{Schur}, we know that  $\widetilde{m}(z):=(\widetilde{m}_i(z))_{i\in J^c}$ satisfies the following perturbed QVE:
\[ -\frac{1}{\widetilde{m}_i(z)} = z+ \sum_{j}^{(J)} X^2_{ij} \widetilde{m}_j(z) + d_i(z), \ i \in J^c\]
where $d_i(z) =Y_i(z) + \sum_{j \in J^c} X_{ij}^2 (G_{jj}^{(iJ)}-G_{jj}^{(J)})$ and $Y_i(z)$ is defined in \eqref{defYi}. On the other hand, let $m(z)=(m_i(z))_{i\in J^c}$  be the unique solution in $\mathbb{H}^{J^c}$ of the unperturbed QVE:
\begin{equation} \label{QDVunpert} - \frac{1}{m_i(z)} = z+ \sum_{j\in J^c} X_{ij}^2 m_j(z), \ i\in J^c.\end{equation}
By definition, $\widetilde{g}(z) =\frac{1}{|J^c|} \sum_{i\in J^c} \widetilde{m}_i(z)$ and $g(z) = \frac{1}{|J^c| }\sum_{i\in J^c} m_i(z)$.  We now specify the stability estimate of Lemma \ref{stability} in our setup where $\mathfrak{X} = J^c$ and $\pi$ is the uniform probability measure on $J^c$. Equation \eqref{QDVunpert} is then the QVE associated to the operator $S: \CC^{J^c}\to \CC^{J^c}$ defined by $(S w)_i =\sum_{j\in J^c} X_{ij}^2 w_j$ for any  $w\in \CC^{J^c}$ and $i\in J^c$. One can easily check that the operator norm of $S$ with respect to the sup norm on $\CC^{J^c}$ is $\|S\|_\infty = \max_{i\in J^c} \|X_i^{(J)}\|^2 \leq C^2$, by definition of $J$ (see \eqref{defJC}). Therefore, we have in particular that $\mathcal{D}_C \subset \{\zeta \in \mathbb{H} : \Im \zeta \geq [\kappa (\|S\|_\infty\vee 1)^2]\vee |\Re \zeta|\}$. 
It follows by Lemma \ref{stability} that 
\begin{equation} \label{controlm} |\widetilde{g}(z)-g(z)|\leq \|\widetilde{m}(z)-m(z)\|_{L^2} \leq \kappa C^2 \|d(z)\|_{L^2}, \end{equation}
where $\| \ \|_{L^2}$ is the $L^2$ norm with respect to the uniform probability measure on $J^c$. 
We will show that for any $\delta>0$,  if $\#\{i\in J^c : |Y_i(z)|\geq \delta^2\}\leq \delta^2 n$ and $|J|\leq n/2$, then
\begin{equation} \label{claimstab}\|d(z)\|_{L^2} \leq 5C^2\delta + R^2/np.\end{equation}
Indeed, observe on the one hand that for any $i \in J^c$,
\[ |Y_i(z)| =\big|\langle X_i^{(J)},G^{(iJ)}X_i^{(J)}\rangle - \sum_{k}^{(J)}X_{ik}^2 G_{kk}^{(iJ)} \big|\leq 2\|X_i\|^2 \|G^{(iJ)}\| \leq 2C^2,\]
where we used that $\|G^{(iJ)}\|\leq 1/\Im z \leq 1$. Thus, on the event where $\#\{i\in J^c : |Y_i|\geq \delta^2\}\leq \delta^2 n$ and $|J|\leq n/2$, we have
 \begin{equation} \label{boundL2Y} \frac{1}{|J^c|} \sum_{i \in J^c} |Y_i(z)|^2 \leq \frac{1}{|J^c|}(2C^2 \delta^2 n + \delta^2|J^c|)\leq 4C^2 \delta^2+\delta^2\leq 5C^2\delta^2.\end{equation}
On the other hand, denote for any $i \in J^c$ by $\widetilde{Y}_i(z):= \sum_{j \in J^c} X_{ij}^2 (G_{jj}^{(iJ)}-G_{jj}^{(J)})$. Using the resolvent identity \eqref{resolvid1} and the fact that $|X_{ij}|\leq R/\sqrt{np}$, we get
\begin{equation} \label{ineqYtilde}  |\widetilde{Y}_i(z)|\leq \sum_j^{(J)} X_{ij}^2 \frac{|G_{ij}^{(J)}|^2}{|G_{ii}^{(J)}|}\leq \frac{R^2}{np} \sum_{j}^{(J)} \frac{|G_{ij}^{(J)}|^2}{|G_{ii}^{(J)}|}.\end{equation}
By Ward's identity $\sum_{j\in J^c} |G_{ij}^{(J)}|^2 = \Im G_{ii}^{(J)}/\Im z$. Therefore, we deduce from \eqref{ineqYtilde} and $\Im z \geq 1$ that $|\widetilde{Y}_i(z)|\leq R^2/np$, for any $i\in J^c$. Together with \eqref{boundL2Y}, this implies that $\|d(z)\|_{L^2}\leq 5C\delta^2+ R^2/np$, and therefore proves the claim \eqref{claimstab}. Putting together \eqref{controlm} and \eqref{claimstab}, it follows that for any $\delta>0$, and $n$ large enough
\[ \PP\big( |\widetilde{g}(z)-g(z)|>6\kappa C^4\delta, |J|\leq n/2\big) \leq \PP\big( \#\{i\in J^c : |Y_i(z)|\geq \delta^2\}> \delta^2 n\big),\]
Invoking Proposition \ref{negli}, this ends the proof.
\end{proof}

\section{Complexity of the kernels $W_n$ for the cut norm}
Equipped with Proposition \ref{equivexpomain1}, we are now reduced to derive a large deviations principle for $(\upsilon_{W_n})_{n\in \NN}$, where $W_n$ is the kernel defined in \eqref{defWn}. Our general strategy is to prove a LDP for $(W_n)_{n\in \NN}$ with respect to a well-chosen topology, and then to contract it to obtain a LDP for $(\upsilon_{W_n})_{n\in \NN}$.
To successfully carry out the contraction principle  
(see \cite[Theorem 4.2.1]{DZ})
 requires two ingredients: firstly, the chosen topology should be coarse enough so that one is able to derive a large deviations principle with a good rate function, and secondly, the topology has to be fine enough so that it renders continuous our function of interest.

In this section, we are concerned with the first point, and prove a key preliminary result. As we will see, the topology induced by the so-called cut norm turns out to be the right one for studying the large deviations of $(W_n)_{n\in \NN}$. Define the cut-norm $\|\ \|_\Box$ on $L^1([0,1]^2)$ as 
\[ \| W\|_\Box = \sup_{S,T \subset [0,1]} \Big|\int_{S\times T} W(x,y) dx dy\Big|, \ W\in L^1([0,1]^2),\]
where the supremum runs over Borel measurable subsets of $[0,1]$. Equivalently, one can use the following functional representation of the cut norm, which will be at times beneficial (see \cite[Lemma 8.10]{Lovasz}), 
\begin{equation} \label{eqcutnorm}\|W\|_\Box = \sup_{f,g :[0,1]\to [0,1]} \Big| \int f(x)g(y) W(x,y) dx dy \Big|, \ W\in L^1([0,1]^2),\end{equation}
where the supremum runs over Borel measurable functions. Recall $\mathcal{W}$ the subset of $L^1([0,1]^2)$ consisting of non negative symmetric functions. Denote by $d_\Box$ the distance induced by $\| \ \|_\Box$ on $\mathcal{W}$, and by $B_{\Box}(W,\delta)$ the closed ball of radius $\delta>0$ centered at $W\in \mathcal{W}$.  

Aiming at understanding the large deviations of $(W_n)_{n\in \NN}$ for the metric induced by the cut norm, we will first show that with overwhelming probability, $W_n$ lives in a subset of $\mathcal{W}$ that can be covered by at most $e^{O(n \log n)}$ balls for the distance $d_\Box$, and thus of negligible $d_\Box$-metric entropy compared to our large deviation speed  $n^2p$ when $np \gg \log n$. This fact will be crucial to derive a LDP upper bound for $(W_n)_{n\in \NN}$ in $(\mathcal{W},d_\Box)$, as it will allow us to reduce the problem to only compute ball probabilities.

More precisely, let $\Pi_m$ be the set of partitions of $(0,1]$ consisting of at most $m+1$ intervals, that is, partitions of the form $\{(0,a_1],(a_1,a_2],\ldots,(a_k,1]\}$, where $0\leq a_1\leq a_2\leq \ldots\leq a_k\leq 1$ with $k\leq m$, and define $\mathcal{K}_{m,r}$ as, 
\begin{equation} \label{defKmr} \mathcal{K}_{m,r} := \big\{\sum_{P,Q\in \mathcal{P}} a_{P,Q} \Car_{P\times Q} \in \mathcal{W} : 0\leq a_{P,Q}\leq r, a_{P,Q}=a_{Q,P}, \mathcal{P} \in \Pi_m\big\}.\end{equation}  
Further, denote by ${S}_n$ the set of permutations of $[n]$ and by $S_{[0,1]}$ the set of Borel measurable bijections of $[0,1]$ preserving the Lebesgue measure. A permutation $\sigma$ of $[n]$ naturally defines an element $\phi_\sigma$ of $S_{[0,1]}$ affine on each interval $I_{in}:=(\frac{i-1}{n},\frac{i}{n}]$, sending $I_{in}$ to $I_{\sigma(i)n}$ for any $i \in [n]$, and such that $\phi_\sigma(0)=0$. Moreover, $S_{[0,1]}$ acts on $\mathcal{W}$ in the following way. If $\phi \in S_{[0,1]}$ and $W\in\mathcal{W}$, we can define the relabelled kernel $W^\phi \in \mathcal{W}$ as
\begin{equation} \label{act} (W^\phi)(x,y) := W(\phi(x),\phi(y)), \ (x,y) \in [0,1]^2.\end{equation}
With a slight abuse of notation, we will write $W^\sigma$, where $W\in \mathcal{W}$ and $\sigma \in S_n$, to denote the relabelled kernel $W^{\phi_\sigma}$, where $\phi_\sigma$ is the measure-preserving bijection of $[0,1]$ associated to $\sigma$. 
With this notation, we have the following result, which can be seen as the counterpart of \cite[Lemma 2.4]{CV} in the sparse case.

\begin{Pro}\label{expopresqtens}
There exists $m=m(r,\veps)$ for any $r,\veps>0$ such that
 \[ \lim_{r\to +\infty} \limsup_{n\to+\infty} \frac{1}{n^2p} \log \PP\big( \inf_{\sigma \in {S}_n} d_\Box (W_n^\sigma , \mathcal{K}_{m,r}) >\veps\big)=-\infty,\]
for any $\veps>0$. 
\end{Pro}

One can check that $\mathcal{K}_{m,r}$ is compact for the $L^1$ norm, and therefore also for the cut norm since $\| \ \|_\Box\leq \|\ \|_1$. In particular, $\mathcal{K}_{m,r}$ is precompact for the distance induced by cut-norm. Since $\log | {S}_n| \sim_{+\infty} n \log n$, Proposition \ref{expopresqtens} entails that $W_n$ is indeed in a subset of $\mathcal{W}$ of $d_\Box$-metric entropy at most $n \log n$, which will be a key ingredient in the proof of the large deviations of $W_n$ for the topology induced by $d_\Box$.

To prove Proposition \ref{expopresqtens} will require a few intermediate steps. Just as in the dense case where $p$ is independent of $n$, the proof relies on the so-called Regularity Lemma. Since the kernels we are considering are unbounded given that $p\ll 1$, we cannot use the usual Regularity Lemma for graphons, instead, we will work with a version of the Regularity Lemma  proven by Borgs, Chayes, Cohn and Zhao \cite{BCCZ} for a class of kernels they called {\em upper regular}. A central object of the Regularity Lemma and of the definition of upper regular kernels is the {\em stepped operator} associated to a partition. Define for any partition $\mathcal{P}$ of $(0,1]$ into Borel subsets $P_1,\ldots,P_k$ and a kernel $W\in \mathcal{W}$, the $\mathcal{P}$-stepped kernel $W_{\mathcal{P}}$ as
\begin{equation} \label{defstepped} (W_\mathcal{P})(x,y) = \frac{1}{\lambda(P_i)\lambda(P_j)} \int_{P_i \times P_j} W(s,t) ds dt, \quad (x,y) \in P_i\times P_j, \ i,j\in [k].\end{equation}
With a slight abuse of notation, for a partition $\mathcal{Q}$ of $[n]$ and a kernel $W$ we will write $W_{\mathcal{Q}}$ to denote the stepped function associated to the partition of $(0,1]$, $\mathcal{P} = \{ \{x\in (0,1] : \lceil nx \rceil \in Q\} : Q\in \mathcal{Q}\}$. Recall the definition of kernels associated to edgeweigthed graphs (see Definition \ref{WG}). With this convention, observe that for any edgeweighted graph $\mathcal{G}$ with vertex set $[n]$ and $\mathcal{P} =\{P_1,\ldots,P_k\}$ a partition of $[n]$ the $\mathcal{P}$-stepped kernel $(W^\mathcal{G})_{\mathcal{P}}$ is described in terms of the normalised edges densities as
\begin{equation} \label{WPG} (W^\mathcal{G})_{\mathcal{P}}(x,y) = d_p(P_i,P_j), \quad (\lceil nx\rceil,\lceil ny\rceil) \in P_i\times P_j,\  i,j\in [k]\end{equation} 
where $d_p(P_i,P_j) = \frac{e_\mathcal{G}(P_i,P_j)}{p|P_i||P_j|}$,  $e_\mathcal{G}(P_i,P_j) = \sum_{u\in P_i, v\in P_j} \beta_{uv}$ and $\beta$ is the adjacency matrix  of $\mathcal{G}$. With this notation, we are now ready to give the definitions of upper regular kernels and upper regular edgeweighted graphs.

\begin{Def}[{\cite[Definition C.2]{BCCZ}}]
\label{upperregular}
Let $\eta>0$ and $K :(0,+\infty) \to (0,+\infty)$ be a Borel measurable function. A kernel $W\in \mathcal{W}$ is said to be {\em upper $(\eta,K)$-regular} if for any partition $\mathcal{P}$ of $(0,1]$ into parts of measure at least $\eta$ and any $\veps>0$,
\begin{equation} \label{upperreg} \| W_\mathcal{P} \Car_{W_\mathcal{P}>K(\veps)} \|_1 \leq \veps.\end{equation}
An edgeweighted graph $\mathcal{G}$ with vertex set $[n]$ is said to be {\em upper $(\eta,K)$-regular} if $W^\mathcal{G}$ is upper $(\eta,K)$-regular, except that only partitions which corresponds to partitions of $[n]$ into parts of size at least $\eta n$ are considered.
\end{Def}

\begin{Rem}
Note that if for some $W\in \mathcal{W}$, $\eta>0$ and $K:(0,1)\to (0,+\infty)$, we have for any $\veps \in (0,1)$, and any partition $\mathcal{P}$ of $(0,1]$ into parts of measure at least $\eta$ that $\|W_\mathcal{P}\Car_{W_\mathcal{P} >K(\veps)}\|_1\leq \veps$, then clearly $W$ is upper $(\eta, \widetilde{K})$- regular with $\widetilde{K}(\veps) := K(\veps\wedge 1)$, $\veps>0$. Thus, it is actually enough to consider in the definition of upper regular kernels functions $K$ defined on $(0,1)$. 
\end{Rem}

For such upper regular kernels, Borgs, Chayes, Cohn and Zhao obtained the following Regularity Lemma. \newpage
\begin{Lem}[{\cite[Theorem C.11]{BCCZ}}]\label{regularity}
Let $K :(0,+\infty)\to (0,+\infty)$ and $\veps>0$. There exists constants $N=N(K,\veps)$ and $\eta_0=\eta_0(K,\veps)$ such that the following holds for any $\eta \leq \eta_0$: for any upper $(K,\eta)$-regular kernel $W$, there exists a partition $\mathcal{P}$ of $[0,1]$ into at most $4^N$ parts each having measure at least $\eta$ such that 
\[ \|W-W_\mathcal{P}\|_\Box \leq \veps.\]
In addition, for any $\eta\leq \eta_0$ and any upper $(K,\eta)$-regular edgeweighted graph with vertex set $[n]$, the partition  $\mathcal{P}$ can be taken to correspond to a partition of $[n]$. 
\end{Lem}

%

%
%
%

In the following proposition, we check that there exists a function $K_\alpha$ such that for any function $\eta_\alpha$, the edgeweighted graph $\mathcal{G}_n$ is upper $(\eta_\alpha,K_\alpha)$-regular with overwhelming probability at the exponential scale $n^2p$ when $\alpha\to +\infty$. 
\begin{Pro}\label{hypoH} Define for any $\alpha\geq 1$, the function $K_\alpha : (0,1)\to (0,+\infty)$ such that $\veps \frac{h_L(K_\alpha(\veps))}{K_\alpha(\veps)} = \alpha$ for any $\veps\in (0,1)$, where $h_L$ is as in \eqref{defL}. Then, for any positive function $\alpha \in [1,+\infty)\mapsto \eta_\alpha$,
\begin{equation} \label{claimupper}\lim_{\alpha\to +\infty}\limsup_{n\to+\infty} \frac{1}{n^2p} \log \PP\big( \mathcal{G}_n \text{ is not upper } (\eta_\alpha ,K_\alpha)\text{-regular}\big) = -\infty.\end{equation}
\end{Pro}
Before proving Proposition \ref{hypoH}, we collect some properties on the function $h_L$ which will be useful as well in the sequel. 

\begin{Lem}The following statements hold:
\label{prophL}
\begin{enumerate}
\item $h_L$ is convex and $h_L(u)=+\infty$ for any $u<0$.
\item $h_L$ is continuously differentiable on $(0,+\infty)$ and $h_L'$ is the inverse of $L'$. 
\item $h_L(0)=1$ and $h_L$ vanishes only at $1$.
\item $h_L(u)\sim_{+\infty} \frac{u}{R^2} \log \frac{u}{R^2}$, where $R$ is the essential supremum of $A_{12}$.

\end{enumerate}
\end{Lem}

\begin{proof}
(1). Recall that $L$ is defined by $L(\theta)=\EE(e^{\theta A_{12}^2})-1$, $\theta\in \RR$, and $h_L$ that is the convex conjugate of $L$. Therefore $h_L$ is by definition convex. Since $L(\theta) \to- 1$ as $\theta \to -\infty$, we deduce that $h_L(u) = +\infty$ for any $u<0$.\\
 (2). Since $A_{12}$ is bounded, $L$ is clearly finite and differentiable on $\RR$ with $L'(\theta) = \EE(A_{12}^2 e^{\theta A_{12}^2})$ for any $\theta \in \RR$. Therefore $L'$ is an increasing homeomorphism from $\RR$ to $(0,+\infty)$. Since $h_L$ is the conjugate of $L$, their subdifferential maps are inverse from each in the sense of multi-valued functions (see \cite[Theorem 23.5]{Rockafellar}). This implies that for any $u>0$, the subdifferential  $\partial h_L(u)$ of $h_L$ at $u$ is reduced to a singleton. As a result, $h_L$ is finite, differentiable on $(0,+\infty)$ and $h_L'$ is the inverse of $L'$ on $(0,+\infty)$. Moreover, since $L'$ is a homeomorphism, $h_L'$ is continuous. \\
(3). We have $h_L(0) =- \inf_\theta L(\theta) = - \lim_{\theta \to -\infty} L(\theta)=1$. Besides, we saw that $L'$ is increasing and $L'(0)=1$ as $\EE(A_{12}^2)=1$. Since $h_L' : (0,+\infty) \to \RR$ is the inverse of $L'$, it follows that $h_L'(u)<0$ for $u\in(0,1)$, $h_L'(1)=0$ and $h_L'(u)>0$ for $u>1$. As a result, $h_L$ vanishes only at $1$.\\
(4). First, note that since $A_{12}^2\leq R^2$ a.s., we have $L(\theta) \leq \phi(\theta R^2)$ for any $\theta \geq 0$, where $\phi(\zeta):=e^\zeta-1$ for any $\zeta\in \RR$. By duality, this entails that $h_L(u) \geq h(u/R^2)$ for any $u\geq 1$, with $h(v)=v \log v-v+1$ for $v\geq 1$. Since $h(v)\sim_{+\infty} v \log v $, this proves the lower bound. To show the upper bound, define $\gamma_\zeta:=\frac{1}{\zeta} \log \EE (e^{\zeta A_{12}^2})$ for any $\zeta >0$. By Hölder's inequality,  $\zeta \mapsto \gamma_\zeta$ is increasing and $\gamma_\zeta \to R^2$ as $\zeta\to +\infty$. Moreover, again by Hölder's inequality $L(\theta) \geq \phi(\gamma_\zeta\theta)$ for any $\theta \geq \zeta$. Fix $\zeta >0$. Using that $L'$ is increasing we can write for any $u \geq L'(\zeta) \vee R^2 e^{R^2\zeta}$, 
\[ h_L(u) = \sup_{\theta \geq \zeta} \{ \theta u- L(\theta)\} \leq  \sup_{\theta \geq \zeta} \{ \theta u- \phi(\gamma_\zeta \theta)\} = h\Big( \frac{u}{\gamma_\zeta}\Big).\]
 Since $\gamma_\zeta \to R^2$ as $\zeta \to +\infty$, this ends the proof of the asymptotic equivalent. 

\end{proof}

We are now ready to prove Proposition \ref{hypoH}.
\begin{proof}[Proof of Proposition \ref{hypoH}] First, we prove that $K_\alpha$ is well-defined and decreasing for any $\alpha\geq 1$. Define $\psi : u\mapsto h_L(u)/u$ for any $u>1$. By Lemma \ref{prophL},  $h_L$ is differentiable on $(0,+\infty)$, convex and $h_L(1)=0$. One can check that this entails that $\psi$ is continuous and increasing on $[1,+\infty)$.   Moreover, again by Lemma \ref{prophL} (3,4), we know that $\psi(u) \to +\infty$ as $u\to+\infty$ and that $\psi(1)=0$. Thus, $\psi : [1,+\infty) \to [0,+\infty)$ is an increasing homeomorphism. As a result,  for any $\alpha\geq 1$, $K_\alpha$ is uniquely defined and decreasing.

To prove  \eqref{claimupper}, we claim that it is enough to show that 
 \begin{equation} \label{claimunion} \lim_{\alpha \to+\infty}\limsup_{n\to+\infty} \frac{1}{n^2p} \max_{\mathcal{P}\in \mathscr{P}_{\alpha}}\max_{\veps \in \mathcal{N}} \log \PP\big(\| (W_n)_{\mathcal{P}} \Car_{(W_n)_{\mathcal{P}}>K_\alpha(2\veps)} \|_{L^1}  >\veps\big) =-\infty,\end{equation}
 where $\mathscr{P}_\alpha$ is the set of partitions of $[n]$ of size $k\leq 1/\eta_\alpha$ and $\mathcal{N} = \{2^{-m} : 0\leq m \leq m_n\}$ with $m_n = \lceil \log_2(\psi(R^2/p)/\alpha)\rceil+1$. Indeed, by Lemma \ref{prophL} (4), we know that $\psi(u) \sim_{+\infty} \log(u/R^2)/R^2$ so that we have $\psi^{-1}(t) \sim_{+\infty} R^2 e^{R^2 t}$. It follows that   $|\mathcal{N}| = O(1/p)$, and besides $|\mathscr{P}_\eta| \leq (1/\eta_\alpha)^n$. Therefore, once the above estimate \eqref{claimunion} is proven, we can use a union bound and the fact that $np \gg 1$ to obtain that 
 \[ \lim_{\alpha \to +\infty} \limsup_{n\to+\infty} \frac{1}{n^2p} \log \PP( \mathcal{E}_{n,\alpha}^c) = -\infty,\]
 where $\mathcal{E}_{n,\alpha}$ is the event where for any partition $\mathcal{P}$ of $[n]$ of size $k\leq 1/\eta_\alpha$ and any $\veps\in \mathcal{N}$, $\| (W_n)_{\mathcal{P}} \Car_{(W_n)_{\mathcal{P}}>K_\alpha(2\veps)} \|_{L^1} \leq \veps$. It now remains to check that on $\mathcal{E}_{n,\alpha}$, $\mathcal{G}_n$ is $(\eta_\alpha,K_\alpha)$-upper regular.  Assume $\mathcal{E}_{n,\alpha}$ occurs and let $\mathcal{P}$ be a partition of $[n]$ into parts of size at least $\eta n$ and $\veps>0$. Let  $\veps \leq 2^{-m_n}$. By definition of $K_\alpha$ and $m_n$, we have $K_\alpha(2\veps) \leq K_\alpha(2^{-m_n+1})  \leq R^2/p$. Since $\| (W_n)_{\mathcal{P}}\|_{L^\infty} \leq R^2/p$, we deduce that $\| (W_n)_{\mathcal{P}} \Car_{(W_n)_{\mathcal{P}}>K_\alpha(2\veps)}\|_{L^1}=0$ so that \eqref{upperreg} trivially holds. Now assume that $\veps> 2^{-m_n}$. There exists $0\leq m \leq m_n$ such that $2^{-m}<  \veps \leq 2^{-m+1}$. Since $\mathcal{P}$ is a partition into parts of size at least $\eta_\alpha n$, the number of parts $k$ has to be less than $1/\eta_\alpha$. Therefore, using in addition the fact that $K_\alpha$ is decreasing and that $\mathcal{E}_{n,\alpha}$ occurs, we get that $\| (W_n)_{\mathcal{P}} \Car_{(W_n)_{\mathcal{P}}>K_\alpha(2\veps)}\|_{L^1} \leq \| (W_n)_{\mathcal{P}} \Car_{(W_n)_{\mathcal{P}}>K_\alpha(2^{-m+1})}\|_{L^1} \leq 2^{-m} <\veps$. Thus, we have shown that on the event $\mathcal{E}_{n,\alpha}$, the edgeweighted graph $\mathcal{G}_n$ is upper $(\eta_\alpha,K_\alpha)$-regular.

We are now reduced to prove \eqref{claimunion}. To this end, fix $\mathcal{P}$ a partition of $[n]$ of size $k\geq 1$. Since $W_n=W^{\mathcal{G}_n}$, we have using \eqref{WPG} that for any $\veps>0$,
\begin{equation} \label{form} \| (W_n)_{\mathcal{P}} \Car_{(W_n)_{\mathcal{P}}>K_\alpha(2\veps)} \|_{L^1} = \frac{1}{n^2}\sum_{1\leq i,j\leq k}d_p(P_i,P_j) \Car_{d_p(P_i,P_j)>K_\alpha(2\veps)} |P_i| |P_j|,\end{equation}
where $d_p(P_i,P_j) = e_{\mathcal{G}_n}(P_i,P_j)/ [p|P_i||P_j|]$ and $e_{\mathcal{G}_n}(P_i,P_j) = \sum_{u\in P_i, v\in P_j} A_{uv}^2\xi_{uv}$. Define for any $r>0$, 
\[ D_{ij}^{(r)}  = d_p(P_i,P_j) \Car_{d_p(P_i,P_j)>r} |P_i||P_j|,\]
and denote by $\Lambda_{ij}^{(r)}$ its log-Laplace transform for any $1\leq i,j\leq k$. We will show that for any $\theta$ such that $L(2\theta) \leq r$, 
\begin{equation} \label{Laplacetrunc}  \Lambda_{ij}^{(r)}(2\theta p) \leq -p|P_i||P_j| (h_L(r) - 2r\theta), \ i,j \in [k], i\neq j,\end{equation}
\begin{equation} \label{Laplacetruncdiag}  \Lambda_{ii}^{(r)}(\theta p) \leq -p \frac{|P_i|^2}{2} (h_L(r) - 2r\theta), \ i\in [k].\end{equation}
For any $\theta' \geq 0$ and $i,j\in [k]$ we have
\begin{align}
\Lambda_{ij}^{(r)}(\theta' p) &= \log\EE( e^{\theta' p  D_{ij}^{(r)}}) =\log \EE(e^{\theta' e_{\mathcal{G}_n}(P_i,P_j) \Car_{d_p(P_i,P_j)>r}})\nonumber \\
&\leq \log \big(1+\EE(\Car_{d_p(P_i,P_j)>r} e^{\theta' e_{\mathcal{G}_n}(P_i,P_j)})\big) \nonumber\\
&\leq \EE(\Car_{d_p(P_i,P_j)>r}e^{\theta' e_{\mathcal{G}_n}(P_i,P_j)}). \label{boundLaplace}\end{align}
 Integrating first on $\xi_{k\ell}$ and using the concavity of the log, we get
\begin{equation} \label{boundLambdap} \log \EE( e^{\theta' \xi_{k\ell}A_{k\ell}^2})  = \log (p \EE (e^{\theta' A_{k\ell}^2}) + 1-p) \leq p L(\theta'),\end{equation}
 for any $\theta'\in \RR$ and $k,\ell\in [n]$. By independence, it follows that for any $i,j \in [k]$, $i\neq j$,
\begin{equation} \label{transfoL} \log \EE\big( e^{\theta' e_{\mathcal{G}_n}(P_i,P_j)}\big)  \leq p| P_i| |P_j| L(\theta'),  \ \log \EE\big( e^{\theta' e_{\mathcal{G}_n}(P_i,P_i)}\big)  \leq p\frac{| P_i|^2}{2}  L(2\theta').\end{equation}
Let $i\neq j$. By Chernoff inequality, we deduce that for any $\theta,\zeta\geq 0$, 
\begin{align*} \log\EE(\Car_{d_p(P_i,P_j)>r} e^{2\theta e_{\mathcal{G}_n}(P_i,P_j)})& \leq  -p|P_i||P_j|r\zeta +\log \EE (e^{(2\theta+\zeta)e_{\mathcal{G}_n}(P_i,P_j)})\\
& \leq -p|P_i||P_j|\big(r\zeta - L(2\theta+\zeta)\big).
\end{align*}
Optimizing on $\zeta \geq 0$, we obtain for any $\theta\geq0$,
\[  \log\EE(\Car_{d_p(P_i,P_j)>r} e^{2\theta e_{\mathcal{G}_n}(P_i,P_j)}) \leq -p |P_i||P_j| \big(\sup_{\zeta \geq 2\theta}\{ r\zeta - L(\zeta)\} -2r\theta\big).\]
Taking $\theta$ such that $L'(2\theta) \leq r$, we get that $\sup_{\zeta \geq 2\theta}\{ r\zeta - L(\zeta)\}   =h_L(r)$, where $h_L$ is defined in \eqref{defL}, which gives \eqref{Laplacetrunc}.
One obtains similarly  using \eqref{transfoL} the estimate \eqref{Laplacetruncdiag}. 

Now, choosing $2\theta_0 = h_L(r)/r$, we find by \eqref{Laplacetrunc}  and \eqref{Laplacetruncdiag} that $\Lambda_{ij}^{(r)}(2p\theta_0)\leq 1$ for any $i\neq j$ and $\Lambda_{ii}^{(r)}(p\theta_0 )\leq 1$ for any $i\in [k]$.
Using Chernoff inequality and the independence of the variables $(D_{ij}^{(r)})_{i\leq j}$, we deduce that 
\begin{align*}
  \log \PP\big( \sum_{1\leq i,j \leq k}D^{(r)}_{ij}>\veps n^2\big) &=\log \PP\big( \sum_{i=1}^k D^{(r)}_{ii} + 2\sum_{i<j} D^{(r)}_{ij}>\veps n^2\big)\\
&\leq - \veps\theta_0 n^2p + \sum_{i=1}^k \Lambda^{(r)}_{ii}(p\theta_0) +\sum_{i<j} \Lambda^{(r)}_{ij}(2p\theta_0)\\
& \leq -\veps\theta_0 n^2p + k^2.
\end{align*}
Taking $r = K_\alpha(2\veps)$, we have by definition of $K_\alpha$ that $\veps \theta_0 = \frac{\veps}{2} \frac{h_L(K_\alpha(2\veps))}{K_\alpha(2\veps)} =  \frac{\alpha}{4}$. Coming back to \eqref{form}, we have show that for any $\veps>0$ and partition $\mathcal{P}$ of $[n]$ of size $k$, 
\[ \PP\big( \|(W_n)_{\mathcal{P}} \Car_{(W_n)_{\mathcal{P}}>K_\alpha(2\veps)} \|_{L^1} > \veps\big) \leq -\frac{\alpha}{4} n^2p +k^2,\] 
which implies \eqref{claimunion}.
\end{proof}


We are now ready to give a proof of Proposition \ref{expopresqtens}. 

\begin{proof}[Proof of Proposition \ref{expopresqtens}] Recall the definitions of $K_\alpha$ in Proposition \ref{hypoH} and of $N(K_\alpha,\veps)$, $\eta_0(K_\alpha,\veps)$ for $\veps>0$ from Lemma \ref{regularity}. Fix $\veps>0$ and set $\eta_\alpha := \eta_0(K_\alpha,\veps)$. We will show that if $\mathcal{G}$ is an upper $(\eta_\alpha,K_\alpha)$-regular edgeweighted graph with vertex set $[n]$, then
\begin{equation} \label{claimprecom} \inf_{\sigma \in {{S}}_n} d_\Box((W^\mathcal{G})^\sigma, \mathcal{K}_{m,r}) \leq 2\veps,\end{equation}
with $r = K_\alpha(\veps)$, $m=4^N$ and $N= N(K_\alpha,\veps)$. Using Proposition \ref{hypoH}, this will end the proof.

Now, assuming that $\mathcal{G}$ is upper $(\eta_\alpha,K_\alpha)$-regular edgeweighted graph with vertex set $[n]$, we deduce from Lemma \ref{regularity} that there exists a partition $\mathcal{P}$ of $[n]$ of size at most $m=4^N$ into parts of size at least $\eta_\alpha n$, such that $\|W^\mathcal{G}-(W^\mathcal{G})_\mathcal{P}\|_\Box \leq \veps$. As $\mathcal{G}$ is upper $(\eta_\alpha,K_\alpha)$-regular, we have also that $\|(W^\mathcal{G})_\mathcal{P} -(W^\mathcal{G})_\mathcal{P}\wedge K_\alpha(\veps) \|_{L^1} \leq \veps$. Since $\| \ \|_\Box \leq \| \ \|_{L^1}$, we get by triangular inequality that 
\[ \| W^\mathcal{G} - (W^\mathcal{G})_{\mathcal{P}}\wedge r \|_\Box \leq 2\veps,\]
with $r=K_\alpha(\veps)$.
We claim that there exists $\sigma \in {S}_n$ such that $((W^\mathcal{G})_\mathcal{P}\wedge r)^\sigma \in \mathcal{K}_{m,r}$. Since the cut norm is invariant by the action of ${{S}}_n$ in the sense that $\|W^\sigma\|_\Box = \|W\|_\Box$ for any kernel $W$, this will end the proof of \eqref{claimprecom}. Note that $((W^\mathcal{G})_\mathcal{P}\wedge r)^\sigma = ((W^\mathcal{G})_\mathcal{P})^\sigma\wedge r = ((W^\mathcal{G})^\sigma)_{\sigma^{-1}\mathcal{P}} \wedge r$ for any $\sigma \in {{S}}_n$, where $\sigma^{-1} \mathcal{P}: = \{\sigma^{-1}(P) : P\in \mathcal{P}\}$. Clearly we can find $\sigma \in {{S}}_n$ so that the elements of $\sigma^{-1}\mathcal{P}$ are intervals of integers, and thus the induced partition of $(0,1]$ belongs to $\Pi_m$. For this choice of $\sigma$, we indeed have $((W^\mathcal{G})_\mathcal{P} \wedge r)^\sigma \in \mathcal{K}_{m,r}$, which ends the proof. 

\end{proof}

\section{Large deviation principle for weighted sparse graphs}

The goal of this section is to prove a large deviations principle for the sequence $(W_n)_{n\in \NN}$ in $(\mathcal{W},d_\Box)$. Before stating the result, we define the following  function $H$ on $\mathcal{W}$, 
\begin{equation} \label{defHW} H(W) := \frac{1}{2} \int_{[0,1]^2}h_L(W(x,y)) dx dy \in [0,+\infty], \ W\in \mathcal{W}\end{equation}
where $h_L$ is defined in \eqref{defL}, and we check in the following lemma,  that $H$ is indeed a rate function on $\mathcal{W}$ with respect to $d_\Box$.
\begin{lemma}\label{lciH}
The function $H$, defined in \eqref{defHW}, is lower semi-continuous on $(\mathcal{W}, d_\Box)$.
\end{lemma}
\begin{proof}Let $W\in \mathcal{W}$. 
Since $h_L$ is the conjugate of $L$, we can write as in the proof of \cite[Lemma 2.1]{CV},
\[
H(W) = \frac{1}{2} \int_{[0,1]^2} \sup_{\theta \in \RR} \{\theta W(x,y) -L(\theta)\} dxdy.\] 
Using a truncation argument and the monotone convergence theorem, we deduce that the following equality holds  \begin{equation} \label{represH} H(W)  = \frac{1}{2} \sup_{\theta \in L^\infty([0,1]^2)} \int_{[0,1]^2} \big\{ \theta(x,y) W(x,y) - L(\theta(x,y))\big\} dx dy.\end{equation}
We will show that for any $\tau>0$, if $U \in \mathcal{W}$ and $(U_n)_{n\in \NN}$ is a sequence of $\mathcal{W}$ converging to $W$ in cut norm such that $H(U_n)\leq \tau$ for any $n\in \NN$, then
\begin{equation} \label{lci} H(U) \leq \liminf_{n\to+\infty} H(U_n).\end{equation}
This will indeed imply that $H$ is lower semi-continuous for the cut norm as it entails that the level sets $\{H\leq \tau\}$, $\tau>0$, of $H$ are closed for the cut norm. Let $U\in \mathcal{W}$ and $(U_n)_{n} \in \mathcal{W}^\NN$ such that $\|U-U_n\|_\Box \to 0$ as $n\to +\infty$ and  $H(U_n)\leq \tau$ for any $n\in \NN$. Observe that as $h_L(u)/u \to +\infty$ as $u\to+\infty$ by Lemma \ref{prophL} (4) and $H(U_n)\leq \tau$ for any $n\in \NN$, the sequence $(U_n)_{n\in \NN}$ is uniformly integrable. Now, by the representation \eqref{represH}, it is enough, in order to prove \eqref{lci}, to show that for any $\theta \in L^\infty([0,1]^2)$, 
\begin{equation} \label{cv} \langle \theta, U_n \rangle \underset{n\to+\infty}{\longrightarrow} \langle \theta, U\rangle.\end{equation} 
This is clearly true for any $\theta$ of the form $\Car_{S\times T}$, where $S$ and $T$ are Borel subsets of $[0,1]$ by the very definition of the cut norm. Now, if $\theta = \Car_E$ where $E$ is a Borel subset of $[0,1]^2$, we can find for any $\veps>0$, a Borel subset $F_\veps = \bigcup_{i=1}^{N_\veps} S_i\times T_i$ of $[0,1]^2$, where $S_i$, $T_i$ are intervals, $(S_i\times T_i)_{i=1}^{N_\veps}$ have pairwise null Lebesgue measure intersection, and $\lambda^2(F_\veps \setminus E)\leq \veps$. We can write for any $n\in\NN$,
\[ |\langle \Car_{E},U -U_n\rangle |\leq |\langle \Car_{F_\veps},U-U_n\rangle | + \sup_{m\in \NN} \langle \Car_{F_\veps\setminus A}, U_m\rangle + \langle \Car_{F_\veps\setminus A}, U\rangle.\]
When $n\to+\infty$, $\langle \Car_{F_\veps},U-U_n\rangle$ goes to zero as $\Car_{F_\veps}$ is almost surely equal to a sum of characteristic functions of products of Borel sets. Besides, since $(U_n)_{n\in \NN}$ is uniformly integrable and $U\in L^1([0,1]^2)$, both $\sup_{m\in \NN} \langle \Car_{F_\veps\setminus E}, U_m\rangle $ and $\langle \Car_{F_\veps\setminus A},U\rangle$ go to $0$ as $\veps \to 0$. This shows that \eqref{cv} holds for $\theta = \Car_E$ for any  Borel subset $E$ of $[0,1]^2$, and thus for any simple function. Now if $\theta \in L^\infty([0,1]^2)$, there exists a sequence of simple functions $(\theta_k)_{k\in\NN}$ converging pointwise to $\theta$ such that $|\theta_k|\leq |\theta|$ for any $k\in \NN$. We have for any $n\in\NN$,
\begin{equation} \label{decouptheta} |\langle \theta, U-U_n\rangle | \leq |\langle \theta_k, U-U_n\rangle| + \sup_{m\in \NN} |\langle \theta-\theta_k,U-U_m\rangle|.\end{equation}
For any $\veps>0$, we can write 
\[ \sup_{m\in \NN}  |\langle \theta-\theta_k,U-U_m\rangle| \leq \veps \sup_{m\in \NN}  \|U-U_m\|_{L^1} +2\|\theta\|_{L^\infty} \sup_{m\in \NN} \int_{|\theta-\theta_k|>\veps} (U+U_m)d\lambda.\]
Since $(U_n)_{n\in\NN}$ is uniformly integrable, it is in particular bounded in $L^1([0,1]^2)$. As $\theta_k$ converges pointwise to $\theta$, it converges in particular in probability to $\theta$, so that using the uniform integrability of $(U_n)_{n\in \NN}$,  taking first the limit in the above inequality when $k\to +\infty$ and then when $\veps \to 0$, we obtain that $\sup_{m\in \NN}  |\langle \theta-\theta_k,U-U_m\rangle|$ converges to $0$ as $k\to+\infty$. Coming back to \eqref{decouptheta} and taking first the limit when $n\to+\infty$ and then the limit when $k\to+\infty$, we get finally the claim. 
\end{proof}

We will prove the following LDP for the sequence $(W_n)_{n\in \NN}$. 
\begin{Pro}\label{LDPW}
The sequence $(W_n)_{n\in \NN}$ satisfies a LDP in $(\mathcal{W},d_\Box)$ with speed $n^2p$ and rate function $H$, defined in \eqref{defHW}.

\end{Pro}

Before going into the proof of Proposition \ref{LDPW}, we show the following non asymptotic large deviation estimate, which will be instrumental in the proof of the upper bound. 
\begin{lemma}\label{lemconv}
For any convex subset $\mathcal{U}\subset \mathcal{W}$ such that $\mathcal{U}\cap L^\infty([0,1]^2)$ is closed in $L^\infty([0,1]^2)$ for the weak-* topology,
\[ \log \PP\big( W_n \in \mathcal{U}\big) \leq -n^2 p \inf_{W\in \mathcal{U}} H(W).\]
\end{lemma}

\begin{proof}Denote for any $r>0$ by $\mathcal{W}_r$ the subset of $\mathcal{W}$ consisting of kernels $W\in \mathcal{W}$ such that $W\leq r$ almost surely, and let $\mathcal{U}$ be a subset of $\mathcal{W}$ satisfying the assumptions of the statement.  Observe that as the entries of $A$ are bounded by $R$, we have $0\leq W_n \leq R^2/p$ almost surely. It trivially implies that $ \PP( W_n \in \mathcal{U}) = \PP( W_n \in \mathcal{U}\cap \mathcal{W}_{R^2/p})$. Now,
for any $\theta \in L^\infty([0,1]^2)$, we can write by Chernoff inequality
\begin{equation} \label{Chernoff} \log \PP\big( W_n \in \mathcal{U}\cap \mathcal{W}_{R^2/p}\big) \leq -n^2 p \inf_{W\in \mathcal{U}\cap \mathcal{W}_{R^2/p}} \langle \theta, W\rangle + \log \EE\big(e^{n^2 p \langle \theta, W_n\rangle}\big).\end{equation}
Denote by $\Lambda_p$ the log-Laplace transform of the common law of the variables $(\xi_{ij}A_{ij}^2)_{i<j}$. By \eqref{boundLambdap}, we know that for any $\zeta \in \RR$, $\Lambda_p(\zeta) \leq pL(\zeta)$. Using independence and denoting by $I_{in} ;=(\frac{i-1}{n},\frac{i}{n}]$ for any $i\in [n]$, we get
\begin{align*} \log \EE\big(e^{n^2 p \langle \theta, W_n\rangle}\big) &= \log \EE\Big( e^{2n^2 \sum_{i<j} \xi_{ij}A^2_{ij} \int_{I_{in}\times I_{jn}} \theta(x,y) dx dy}\Big)\\
& = \sum_{i<j} \Lambda_p\Big( 2n^2 \int_{I_{in}\times I_{jn}} \theta(x,y) dx dy\Big)\\
& \leq p\sum_{i<j}L\Big(  2n^2 \int_{I_{in}\times I_{jn}} \theta(x,y) dx dy\Big).
\end{align*}
Since $L$ is convex, it yields by Jensen inequality that 
\[ \log \EE\big(e^{n^2 p \langle \theta, W_n\rangle}\big)  \leq \frac{1}{2} n^2p\int_{[0,1]^2} L(2\theta(x,y)) dx dy.\]
Using this estimate in \eqref{Chernoff} and optimizing on $\theta \in L^\infty([0,1]^2)$ we obtain
\[ \log \PP( W_n \in \mathcal{U}) \leq -n^2 p \sup_{\theta \in L^\infty([0,1]^2)} \inf_{W\in \mathcal{U}\cap \mathcal{W}_{R^2/p}} \big\{ \langle \theta,W\rangle -  \frac{1}{2} \int_{[0,1]^2} L(2\theta(x,y))dx dy\big\}.\]
We are now in position to apply a minimax theorem. The sets $L^\infty([0,1]^2)$ and $\mathcal{U}\cap \mathcal{W}_{R^2/p}\subset L^\infty([0,1]^2)$ are both convex subsets and the functions $\theta \mapsto \langle \theta,W\rangle - \int_{[0,1]^2} L(2\theta(x,y))dx dy$ and $W\mapsto \langle \theta, W\rangle$ are respectively concave and convex. Moreover,  $\mathcal{U}\cap \mathcal{W}_{R^2/p}$ is a weak-* closed  subset of  $(R^2/p) B_{L^\infty([0,1]^2)}$ since $\mathcal{U}\cap L^\infty([0,1]^2)$ is  weak-* closed in $L^\infty([0,1]^2)$ by assumption. As $(R^2/p) B_{L^\infty([0,1]^2)}$ is weak-* compact by Banach-Alaoglu theorem, it follows that $\mathcal{U}\cap \mathcal{W}_{R^2/p}$ is weak-* compact. In addition, for any given $\theta \in L^{\infty}([0,1]^2)\subset L^1([0,1]^2)$, $W \in L^\infty([0,1]^2) \mapsto \langle \theta, W\rangle$ is weak-* continuous. By Ky Fan minimax Theorem \cite[Theorem 4.36]{Clarke}, we deduce that 
\[ \log \PP( W_n \in \mathcal{U}) \leq -n^2 p  \inf_{W\in \mathcal{U}\cap \mathcal{W}_{R^2/p}}\sup_{\theta \in L^\infty([0,1]^2)}  \big\{ \langle \theta,W\rangle -   \frac{1}{2}\int_{[0,1]^2} L(2\theta(x,y))dx dy\big\}.\]
Using \eqref{represH}  and the fact that $\inf_{\mathcal{U}\cap \mathcal{W}_{R^2/p}} H \geq \inf_{\mathcal{U}} H$, this ends the proof.

\end{proof}
Equipped with Lemma \ref{lemconv}, we can now give a proof of Proposition \ref{LDPW}.

\begin{proof}[Proof of Proposition \ref{LDPW}]
\textbf{Upper bound}. Let $F\subset \mathcal{W}$ be a closed subset for the distance induced by the cut norm.
We first make use of Corollary \ref{expopresqtens} to reduce ourselves to case where $F$ is a ball. Recall the definition of $\mathcal{K}_{m,r}$ in \eqref{defKmr}. As observed after the statement of Corollary \ref{expopresqtens}, the set $\mathcal{K}_{m,r}$ is compact for the $L^1$ norm, and as a consequence also for the cut norm since $\| \ \|_\Box \leq \| \ \|_{L^1}$. Thus, there exists $U_1,\ldots, U_N \in \mathcal{W}_r$, with $N$ depending on $r$ and $\veps$, such that $\mathcal{K}_{m,r} \subset \bigcup_{i=1}^{N} B_{\Box}(U_i,\veps)$. This covering, and the invariance of the cut norm by the action of $S_n$ allows us to write that 
\[ \{ \inf_{\sigma \in {S}_n} d_\Box(W_n^\sigma, \mathcal{K}_{m,r})\leq \veps\} \subset  \{W_n \in \bigcup_{\sigma \in {S}_n} \bigcup_{i=1}^N B_{\Box}(U_i^\sigma ,2\veps)\}.\]
Now, if $ F \cap B_{\Box}(U_i^\sigma,2\veps) \neq \emptyset$ for some $i$ and $\sigma$,  then this means that $d_{\Box}(U_i^\sigma,F)\leq 2\veps$. Let $T \subset S_n\times [N]$ be the set of pairs $(\sigma,i)$ such that $F \cap B_{\Box}(U_i^\sigma,2\veps) \neq \emptyset$. Then,
\[ \{ \inf_{\sigma \in {S}_n} d_\Box(W_n^\sigma, \mathcal{K}_{m,r})\leq 2\veps, {W}_n \in F\} \subset  \{W_n \in \bigcup_{(\sigma,i) \in T} B_{\Box}(U_i^\sigma,2\veps)\}.\] Using Corollary \ref{expopresqtens}, a union bound and the fact that $\log | {S}_n |\leq n \log n \ll n^2p$, we deduce that it suffices to show that 
 \begin{equation} \label{claimLD} \lim_{\veps \downarrow 0} \limsup_{n\to+\infty} \frac{1}{n^2p}\sup_{{U} \in V_{\veps}(F)} \log \PP(W_n \in B_{\Box}(U,\veps)) \leq - \inf_F H,\end{equation}
 where $V_{\veps}(F)$ denotes the closed $\veps$-neighborhood of $F$ for the distance $d_\Box$, that is $V_{\veps}(F):=\{W\in \mathcal{W} : d_\Box(W,F) \leq \veps\}$.
Now, note that for any $U\in \mathcal{W}$, $B_{\Box}(U,\veps)$ is a convex set, and $B_\Box(U,\veps) \cap L^\infty([0,1]^2)$ is weak-* closed. By Lemma  \ref{lemconv}, we deduce that for any $U \in \mathcal{W}$, 
\[  \log \PP(W_n \in B_{\Box}(U,\veps)) \leq -n^2 p \inf_{B_{\Box}(U,\veps)} H.\]  
But, $B_{\Box}({U},\veps) \subset V_{2\veps}(F)$ whenever $d_\Box({U},F)\leq \veps$. 
Therefore,
\begin{equation} \label{nonasym}  \sup_{{U} \in V_\veps(F)}\log \PP(W_n \in B_{\Box}(U,\veps)) \leq -n^2 p \inf_{V_{2\veps}(F)} H.\end{equation}
As $H$ is lower semi continuous on $({\mathcal{W}}, d_\Box)$ by Lemma \ref{lciH} and $F$ is a closed set, $\lim_{\veps \downarrow 0} \inf_{V_{2\veps}(F)} H = \inf_F H$. Thus, dividing by $n^2p$ in \eqref{nonasym}, taking the $\limsup$ as $n\to+\infty$ and then letting $\veps \downarrow 0$, we get the claim \eqref{claimLD}, which ends the proof of the upper bound.

\textbf{Lower bound}. We will prove the lower bound by using a classic exponential tilting and Lemma \ref{chgmes}.  Denote by $H(\QQ \mid \PP)$ the relative entropy of a probability measure $\QQ$ with respect to $\PP$, defined as 
\begin{equation} \label{defentropy} H(\QQ \mid \PP) =  \int \frac{d\QQ}{d\PP} \log \frac{d\QQ}{d\PP} d\PP \end{equation}
if $\QQ$ is absolutely continuous with respect to $\PP$ and  where $\frac{d\QQ}{d\PP}$ is the Radon-Nikodym derivative of $\QQ$ with respect to $\PP$, and otherwise defined as $+\infty$.  Now, to prove the large deviations lower bound, it suffices to show that for any $U\in \mathcal{W}$, 
\begin{equation} \label{lb} \lim_{\veps \downarrow 0} \liminf_{n\to+\infty} \frac{1}{n^2p}\log \PP\big(d_\Box(W_n,U) \leq \veps) \geq -H(U).\end{equation}
First, one can easily check that it is enough to prove \eqref{lb} in the case where both $U$ and $1/U$ are bounded. Indeed, denoting by $U_r:= (U\wedge r)\vee r^{-1}$ for any $r\geq 1$ and $U\in \mathcal{W}$, we have $H(U_r) \leq H(U)$ and $(U_r)_{r\in \NN}$ converges in $L^1$ norm to $U$, and thus also for the cut norm. Fix now $U$ such that $U,1/U \in L^\infty([0,1]^2)$. Let $u_{ij}$ denote the average of $U$ on $I_{in}\times I_{jn}$ for any $i\neq j$ and $u_{ii}=0$ for any $i\in [n]$. Set $U_n$ to be the kernel constant on each $I_{in} \times I_{jn}$ and equal to $u_{ij}$
 for any $i,j \in [n]$. It is a classical fact that $(U_n)_{n\in \NN}$ converges to $U$ for the $L^1$ norm and as a result for the cut norm as well. Therefore, it is enough to show that
\begin{equation} \label{lb2} \lim_{\veps \downarrow 0} \liminf_{n\to+\infty} \frac{1}{n^2p}\log \PP\big(d_\Box(W_n,U_n) \leq \veps) \geq -H(U).\end{equation}
Denote by $E_{n,\veps}$ the event $\{ d_\Box(W_n,U_n) \leq \veps\}$ and define $\QQ$ as the probability measure proportional to $\exp\big(\sum_{i<j} \theta_{ij} \xi_{ij} A_{ij}^2\big) d\PP$ where $\theta_{ij} = h_L'(u_{ij})$ for any $i,j\in [n]$. By Lemma \ref{chgmes}, we have
\begin{equation} \label{chgmeseq} \PP(E_{n,\veps}) \geq \QQ(E_{n,\veps}) \exp\Big( -\frac{1}{\QQ(E_{n,\veps})} \big( H(\QQ\mid \PP) + e^{-1}\big)\Big).\end{equation}
We claim that 
\begin{equation} \label{claimentropy} H(\QQ\mid \PP) \leq n^2 p H(U) + O(n^2p^2),\end{equation}
and that the event $E_{n,\veps}$ is typical under the measure $\QQ$, that is, for any $\veps>0$,
\begin{equation} \label{claimevent} \QQ(E_{n,\veps}) \underset{n\to+\infty}{\longrightarrow} 1.\end{equation} 
Once these two statements are proven, the claimed lower bound \eqref{lb2} follows immediately from \eqref{chgmeseq}.
We start by proving \eqref{claimentropy}. Recall $\Lambda_p$, the log-Laplace transform of the joint law of the variables $(\xi_{ij} A_{ij}^2)_{i<j}$. Using independence, we obtain that $H(\QQ\mid \PP) = \sum_{i<j} \big( \theta_{ij}\Lambda'_p(\theta_{ij}) - \Lambda_p(\theta_{ij})\big)$. One can easily check that for any $\theta \in \RR$, 
\begin{equation} \label{boundL} |\Lambda_p(\theta) - pL(\theta)| \leq p^2 a(\theta), \ |\Lambda_p'(\theta)-pL'(\theta)|\leq p^2b(\theta),\end{equation}
where $a$ and $b$ are some continuous functions. Since $U$ and $1/U$ are bounded, the family $(u_{ij})_{i,j\in[n]}$ is uniformly bounded from above and below, and as $\theta_{ij} = h_L'(u_{ij})$ for any $i,j\in[ n]$ and $h_L'$ is continuous by Lemma \ref{prophL}, it follows that  the family  $(\theta_{ij})_{i,j\in [n]}$ is uniformly bounded as well. Thus, by \eqref{boundL} we obtain that
\begin{align} H(\QQ \mid \PP) &= p\sum_{i<j} \big(\theta_{ij} L'(\theta_{ij})-L(\theta_{ij})\big) + O(p^2n^2) \nonumber \\
&= p\sum_{i<j} h_L(u_{ij}) + O(p^2n^2),\label{entropy} \end{align}
where we used the fact that $h_L$ is the conjugate of $L$ and that as $L'$ is the inverse of $h_L'$ by Lemma \ref{prophL}, $L'(\theta_{ij})=u_{ij}$ for any $i,j\in[n]$. By Jensen's inequality we have $h_L(u_{ij}) \leq n^2 \int_{I_{in}\times I_{jn}} h_L(U(x,y))  dx dy$ for any $i<j$. Therefore, \eqref{entropy} implies  that $H(\QQ\mid \PP) \leq n^2p H(U) + O(n^2p^2)$, which ends the proof of \eqref{claimentropy}.

 We now turn our attention to \eqref{claimevent}.  Since both $W_n$ and $U_n$ are constant on each $I_{in}\times I_{jn}$, $i,j\in [n]$ and vanish on $I_{in}\times I_{in}$ for any $i\in [n]$, it follows from \cite[Lemmas E.1, E.2]{Janson} that $d_\Box(W_n,U_n) \leq 6 n^{-2} \max_{S \subset [n]} |\sum_{i\in S,j\in S^c} (\frac{\xi_{ij}}{p}A_{ij}^2-u_{ij})\big|$. (The choice of taking $u_{ii}=0$ for any $i\in[n]$ was motivated by the use this inequality). 
As there are $2^n$ subsets of $[n]$ and $np \gg \log n$, it is enough to show that for $n$ large enough and any $S\subset [n]$, 
\begin{equation} \label{claimconcS} \QQ\Big(\big| \sum_{i\in S,j\in S^c} \big(\frac{\xi_{ij}}{p}A_{ij}^2-u_{ij}\big)\big| > \veps n^2\Big)\leq 2e^{- \frac{n^2p M}{R^2} h( \frac{\veps}{2M})},\end{equation}
where $h$ is defined in Corollary \ref{LDPER} and $M>0$ is some constant depending on $\|U\|_{L^\infty}$ and $\|1/U\|_{L^\infty}$.
Note that $\EE_\QQ(\xi_{ij}A_{ij}^2) = \Lambda_p'(\theta_{ij})$ for any $i,j\in [n]$. Therefore, \eqref{boundL} and the fact that $L'(\theta_{ij})=u_{ij}$ imply that 
\begin{equation} \label{meanQ} \big| \EE_\QQ(\xi_{ij}A_{ij}^2) - p u_{ij}\big| \leq K p^2, \ i,j\in [n],\end{equation}
where $K>0$ depends  on $\|U\|_{L^\infty}$ and $\|1/U\|_{L^\infty}$. Moreover, $\EE_\QQ[(\xi_{ij}A_{ij}^2)^2]\leq R^2 \EE_\QQ(\xi_{ij}A_{ij}^2) \leq R^2(pu_{ij}+Kp^2)\leq M R^2p$ for any $i,j \in [n]$, where $M$ is some constant depending again on $\|U\|_{L^\infty}$ and $\|1/U\|_{L^\infty}$. Since under $\QQ$, $(\xi_{ij}A_{ij}^2)_{i\in S,j\in S^c}$ are independent random variables bounded by $R^2$ and $v:=\sum_{i\in S,j\in S^c} \EE_\QQ[(\xi_{ij}A_{ij}^2)^2] \leq n^2p MR^2$, we deduce by Bennett's inequality (see \cite[Theorem 2.9]{BLM}) that 
\[ \QQ\Big(\big| \sum_{i\in S,j\in S^c} \big(\xi_{ij} A_{ij}^2-\EE_\QQ(\xi_{ij}A_{ij}^2 )\big)\big| > \frac{\veps}{2} n^2p\Big) \leq 2e^{- \frac{n^2p M}{R^2} h( \frac{\veps}{2M})}.\]
Combining this concentration inequality with \eqref{meanQ}, we obtain \eqref{claimconcS}. This finally ends the proof of the lower bound. 
\end{proof}

With Proposition \ref{LDPW}, we have fulfilled the first part of the assumptions of the contraction principle, except for the fact that the rate function $H$ of the LDP of $(W_n)_{n\in \NN}$ has no reason to be a good rate function. To fix this issue, we contract the LDP of $(W_n)_{n\in \NN}$ in $(\mathcal{W},d_\Box)$ to the quotient space $(\widetilde{\mathcal{W}},\delta_\Box)$ of {\em unlabelled kernels}.

More precisely, recall ${S}_{[0,1]}$ the set of Borel measurable bijections of $[0,1]$ preserving Lebesgue measure, and its action on $\mathcal{W}$ described in \eqref{act}. Define the cut metric $\delta_\Box$ by
\[\delta_\Box(W_1,W_2)= \inf_{\phi \in {S}_{[0,1]}} \| W_1-W_2^\phi\|_\Box = \inf_{\psi \in {S}_{[0,1]}}\|W_1^\psi-W_2\|_\Box, \ W_1,W_2\in \mathcal{W},\]
which is a pseudo-metric on $\mathcal{W}$ (see \cite[Lemma 6.5]{Janson}). This implies that the relation $W_1\sim W_2 \Longleftrightarrow \delta_\Box(W_1,W_2)=0$ is an equivalence relation on $\mathcal{W}$,  and as a result we can define the quotient set $\widetilde{\mathcal{W}}$ whose elements are called {\em unlabelled kernels}. The cut metric defines a distance as well on $\widetilde{\mathcal{W}}$ by setting,
$\delta_\Box(\widetilde{W}_1,\widetilde{W}_2) = \delta_\Box(W_1,W_2)$, $W_1,W_2\in \mathcal{W}$, which is compatible with the quotient topology. Now,  the invariance of $H$ with respect to  the action of ${S}_{[0,1]}$ and its lower semi-continuity with respect to the cut norm allows us to define $H$ as well on the quotient space $\widetilde{\mathcal{W}}$ as described in the following lemma. 
\begin{lemma}\label{lciHtilde}
Let for any $\widetilde{W} \in \widetilde{\mathcal{W}}$, $H(\widetilde{W}) : = H(W)$. Then $H$ is well-defined on $\widetilde{\mathcal{W}}$ and lower semi-continuous with respect to $\delta_{\Box}$. 
\end{lemma}
The proof of this lemma is essentially the same as the one of \cite[Lemma 2.1]{CV}, so that we omit it. The interest of moving to the set of unlabelled kernels comes from the fact that tractable compactness criteria are then available. In particular, as proven in \cite[Theorem C.7]{BCCZ}, uniformly integrable subsets of kernels are compact under the cut metric (more generally uniformly upper regular subset of kernels by \cite[Theorem C.13]{BCCZ}). In particular, as we now show, the rate function $H$ has compact level sets on the quotient space $\widetilde{\mathcal{W}}$.

%
%
%

\begin{lemma}\label{goodH}
$H$ is a good rate function on $(\widetilde{\mathcal{W}},\delta_\Box)$. 
\end{lemma}

\begin{proof}
Let $\tau>0$. Since $h_L(u)/u\to +\infty$ as $u\to+\infty$ by Lemma \ref{prophL}, it follows that $\{ H\leq \tau\}$ is a set of uniformly integrable kernels. By \cite[Theorem C.7]{BCCZ}, it follows that $\{H\leq \tau\}$ is relatively compact for the cut metric. Since $\{ H\leq \tau\}$ is in addition closed by Lemma \ref{lciHtilde}, this ends the proof. 
\end{proof}

Finally, using the continuity of the map $W\in \mathcal{W} \mapsto \widetilde{W}\in \widetilde{\mathcal{W}}$ and the contraction principle (see \cite[Theorem 4.2.1]{DZ}), we can deduce from Proposition \ref{LDPW} the following large deviations principle. 

\begin{Pro}
\label{upper}
The sequence $(\widetilde{W}_n)_{n\in\NN}$ satisfies a LDP in $(\widetilde{\mathcal{W}},\delta_\Box)$ with speed $n^2p$ and good rate function $H$ defined in \eqref{defHW}.
\end{Pro}

\section{A counting lemma for trees}
We now come to the last step of our strategy of using the contraction principle to obtain a LDP for $(\upsilon_{W_n})_{n\in \NN}$. Recall Definition \ref{specW} of the QVE measure $\upsilon_W$ of a kernel $W\in \mathcal{W}$ and of the degree truncated kernel $W^{(C)}$, where $C>0$, as 
\begin{equation} \label{defWCrappel} W^{(C)}(x,y) := W(x,y) \Car_{d_W(x)\leq C} \Car_{d_W(y)\leq C}, \ (x,y)\in [0,1]^2,\end{equation}
where $d_W$ is the degree function of $W$, that is $d_W(x):=\int_0^1 W(x,y) dy$, $x\in [0,1]$.
The goal of this section is to prove that the QVE measure is well-defined for any kernel in $\mathcal{W}$ and that the resulting map is continuous for the topology induced by the cut norm.

\begin{Pro}
\label{defcont}
For any $W\in \mathcal{W}$, let $\upsilon_{W} := \lim_{C\to +\infty} \upsilon_{W^{(C)}}$, where $W^{(C)}$ is defined in \eqref{defWC}.
The map ${W} \mapsto \upsilon_{{W}}$ is well-defined and continuous from $({\mathcal{W}},d_{\Box})$ to $\mathcal{P}(\RR)$ endowed with the weak topology.
\end{Pro}

To prove Proposition \ref{defcont}, several steps will be needed. First, recall that we denoted by $\mathcal{X}$ the set of Borel measurable  non negative symmetric  kernels on $[0,1]^2$ with a bounded degree function, where kernels agreeing almost everywhere are identified.  We define further $\mathcal{X}_C$ for any $C>0$ as the subset of $\mathcal{X}$ consisting of kernels $W$ such that $d_W\leq C$ almost surely, and we investigate the question of the continuity of the map $W\in \mathcal{X}_C \mapsto \upsilon_W$. To this end, we will show that the moments of $\upsilon_W$ for $W\in \mathcal{X}$ can actually be expressed in terms of {\em homomorphism densities} of trees. Define for any loopless graph $F=(V(F),E(F))$ and $W\in \mathcal{W}$, the homomorphism density $t(F,W)$ as 
\begin{equation} \label{defhom} t(F,W) := \int_{[0,1]^{V(F)}} \prod_{ij \in E(F)} W(x_i,x_j) \prod_{i\in V(F)} dx_i \in [0,+\infty].\end{equation}
Further, for any rooted loopless graph $(F,o)$, we set
\begin{equation} \label{defhomroot} t_{x_o}(F,W):=\int_{[0,1]^{V(F)\setminus o}} \prod_{ij\in E(F)} W(x_i,x_j) \prod_{i\in V(F)\setminus o} dx_i, \ W\in \mathcal{W}, x_o\in[0,1].\end{equation}
In particular,  with these definitions one has that $\int_0^1 t_{x_o} (F,W)dx_o = t(F,W)$. 
Using \cite[Lemma 2.2]{EM}, we can compute the moments of the QVE measures of kernels in $\mathcal{X}$ as follows. 
\begin{Pro}\label{momentuW}
Let $W\in \mathcal{X}$. Then for any $k\in \NN$, 
\[ \int \tau^{2k} d\upsilon_W(\tau) = \sum_{(F,o)\in \mathcal{T}_{k}} t_{}(F,W), \quad \int \tau^{2k+1} d\upsilon_W(\tau)=0,\]
where $\mathcal{T}_k$ denotes the set of unlabelled rooted planar trees $(F,o)$ with $k$ edges. 
\end{Pro}

\begin{proof}Fix $W\in \mathcal{X}$ and let $(\upsilon_x)_{x\in [0,1]}$ denote the family of probability measures such that for any $z\in \mathbb{H}$,  the family of their Stieltjes transforms $m(z)=(m(z,x))_{x\in [0,1]}$ is the unique solution in $\mathcal{B}^+$ of the QVE \eqref{QVE} associated to $W$. In \cite[Lemma 2.2]{EM}, the moments of $\upsilon_x$ were computed in the case of a finite dimensional QVE, which corresponds to the case where $W$ is a stepped kernel associated to a partition $\mathcal{P}_n:=\{(\frac{i-1}{n},\frac{i}{n}] : i\in[n]\}$, $n\in \NN$. One can check that the exact same proof carries out to the case of general kernels $W\in \mathcal{X}$.  This yields that for any $x\in [0,1]$,
\[ \int \tau^{2k} d\upsilon_x(\tau) = \sum_{(F,o)\in \mathcal{T}_{k}} t_{x_o}(F,W) \Car_{x_o=x}, \quad k\in \NN.\]
Integrating the above inequality gives the claim for the even moments since by definition $\upsilon_W = \int_0^1 \upsilon_x dx$. Besides, by \cite[Theorem 2.1]{AEK}, we know that  $\upsilon_x$ is symmetric for any $x\in [0,1]$. Thus, $\upsilon_W$ is symmetric as well, and as a result its odd moments vanish. 
\end{proof}

In the next lemma, we show that homomorphism densities of trees are continuous for the cut norm, provided a uniform bound on the degree functions of the kernels. First, extend the definitions of homomorphism densities to $\mathcal{W}$-decorated graphs $(F,w)$ and $\mathcal{W}$-decorated rooted graphs $((F,o),w)$ as follows
\[ t(F,w) := \int_{[0,1]^{V(F)}} \prod_{ij \in E(F) } W_{ij}(x_i,x_j) \prod_{i\in V(F)} dx_i \in [0,+\infty],\]
\[ t_{x_o}(F,w) :=  \int_{[0,1]^{V(F)\setminus o}} \prod_{ij \in E(F) } W_{ij}(x_i,x_j) \prod_{i\in V(F)\setminus o} dx_i, \ x_o\in [0,1],\]
where $F=(V(F),E(F))$ is a loopless graph and $w=(W_e)_{e\in E(F)} \in \mathcal{W}^{E(F)}$. By convention, if $F$ is reduced to the vertex $o$, then we set $t_{x_o}(F,w)=1$ for any $x_o\in [0,1]$. 
With this notation, we prove the following counting lemma for trees decorated with kernels in $\mathcal{X}$. 

\begin{lemma}[Counting lemma for decorated trees]\label{countinglem}

Let $(F,w)$ and $(F,w')$ be two $\mathcal{X}$-decorated trees with the same underlying tree $F$, where $w=(W_e)_{e\in E(F)}\in \mathcal{X}^{E(F)}$ and $w'=(W'_e)_{e\in E(F)}\in \mathcal{X}^{E(F)}$. Then,
\[ | t(F,w) - t(F,w')| \leq M^{e(F)-1} \sum_{e\in E(F)} \|W_e-W_e'\|_{{ \Box}},\]
where $M= \max\big(\max_{e\in E(F)}  \|d_{W_e}\|_{L^\infty} , \max_{e\in E(F)} \| d_{W'_e}\|_{L^\infty})$ and $e(F)=|E(F)|$.

\end{lemma}

To prove Lemma \ref{countinglem}, we will need the following a priori bound on the homomorphism densities of $\mathcal{X}$-decorated rooted trees.  
\begin{lemma}\label{boundhomotree}
For any  $\mathcal{X}$-decorated rooted tree $((F,o),w)$, 
\begin{equation} \label{ineqt}  t_{x_o}(F,w) \leq \Big(\max_{e\in E(F)} \|d_{W_e}\|_{L^\infty} \Big)^{e(F)}, \quad x_o \in [0,1],\end{equation}
where $e(F) = |E(F)|$.
\end{lemma}

\begin{proof}
We proceed by induction over the number of edges. If $F$ consists of only its root $o$, then the claim trivially holds since $t_{x_o}(F,w)=1$ by convention. Assume now that the inequality \eqref{ineqt} is true for any $\mathcal{X}$-decorated rooted tree $((F',o'),w')$ such that $e(F')=n$ for a given $n\in \NN$ and let $((F,o),w)$ be a $\mathcal{X}$-decorated rooted tree with $e(F) =n+1$. Since $e(F)\geq1$, there exists a leaf $v\neq o$ in $F$. Let $v'$ denote its unique neighbour in $F$ and set  $F'$ as the tree obtained from $F$ by removing the vertex $v$ and the edge $vv'$. Then, we can write
\[
 t_{x_o}(F,w) = \int_{[0,1]^{V(F)\setminus  o} }\prod_{ij \in E(F')} W_{ij}(x_i,x_j) W_{vv'}(x_v,x_{v'}) \prod_{\ell \in V(F) \setminus o} d{x_\ell}.\]
Since $v\notin F'$ and $o\neq v$, we obtain by integrating first on $x_v$ that
\[  t_{x_o}(F,w)  \leq M  t_{x_o}(F',w'),\]
where $w'=(W_e)_{e\in E(F')}$ and  $M = \max_{e\in E(F)}\|d_{W_e}\|_{L^\infty}$. Together with the induction hypothesis, this ends the proof.
\end{proof}
We are now ready to give a proof of Lemma \ref{countinglem}.
\begin{proof}[Proof of Lemma \ref{countinglem}]
Clearly it suffices to prove the statement when $w$ and $w'$ differ on only one edge, say $uv\in E(F)$. As $F$ is a tree, after removing the edge $uv$, one is left with two disjoint trees $F'$ and $F''$, $F'$ containing $u$ and $F''$ containing $v$. With this notation, we can write
\[ t(F,w) - t(F,w')=  \int_{[0,1]^2} t_{x_u}(F',w') (W_{uv}(x_u,x_v)-W'_{uv}(x_u,x_v)) t_{x_v}(F'',w'')dx_u dx_v,\]
where $w'=(W_e)_{e\in E(F')}$, $w''=(W_e)_{e\in E(F'')}$, $F'$ being rooted at $u$ and $F''$ at $v$.
By Lemma \ref{boundhomotree}, we know that for any $x_u\in [0,1]$, $t_{x_u}(F',w') \leq M^{e(F')}$ and for any $x_v\in [0,1]$, $t_{x_v}(F'', w'') \leq M^{e(F'')}$. 
Thus, using \eqref{eqcutnorm}, we get that $|t(F,w)-t(F,w')|\leq M^{e(F')+e(F'')} \|W_{uv}-W'_{uv}\|_{\Box}$. As $e(F) =e(F)+e(F')+1$, this ends the proof.

\end{proof}

Building on the Counting Lemma \ref{countinglem} and Lemma \ref{momentuW}, we obtain the following continuity result. 
 \begin{Pro}\label{conti}For any $C>0$,
the map $W\mapsto \upsilon_W$ is continuous from $(\mathcal{X}_C,d_\Box)$ to $\mathcal{P}(\RR)$ endowed with the weak topology. 
\end{Pro}

\begin{proof}
Since every probability measure $\upsilon_W$ for any $W\in \mathcal{X}_C$ is compactly supported and symmetric by \cite[Theorem 2.1]{AEK}, it suffices to prove that for any $k\in \NN$, $W\in \mathcal{X}_C\mapsto \int \tau^{2k} d\upsilon_W(\tau)$ is continuous for the cut norm. Using Corollary \ref{momentuW} and Lemma \ref{countinglem}, we deduce that for any $W,W'\in \mathcal{X}_C$, 
\[ \big| \int \tau^{2k} d\upsilon_W(\tau)-\int \tau^{2k} d\upsilon_{W'}(\tau)\big| \leq  k C^{k-1} |\mathcal{T}_k| \|W-W'\|_{\Box},\]
which shows that $W\in \mathcal{X}_C\mapsto \int \tau^{2k} d\upsilon_W(\tau)$ is indeed continuous for the cut norm.
\end{proof}

In order to lift this continuity result on $\mathcal{X}_C$, $C>0$ to the whole space of integrable kernels $\mathcal{W}$, we will need the following version of Hoeffman-Wielandt inequality for QVE measures of kernels. 

\begin{Lem}[Hoeffman-Wielandt inequality]\label{HW}
For any $W,W'\in \mathcal{X}$, 
\[ \mathscr{W}_2(\upsilon_W,\upsilon_{W'}) \leq  \|W-W'\|_{L^1}^{1/2}.\] 
\end{Lem}

The proof of this inequality consists in taking the limit when the dimension goes to infinity in Hoeffman-Wielandt inequality (see \cite[Lemma 2.1.19]{AGZ}) and can be found in the Appendix \ref{appendix}.

\begin{proof}[Proof of Proposition \ref{defcont}]
We first show that $\upsilon_W$ is well-defined for any $W\in \mathcal{W}$ as the limit for the weak topology of $\upsilon_{W^{(C)}}$ as $C\to +\infty$. Denote by $d_{bL}$ the bounded Lipschitz metric on $\mathcal{P}(\RR)$, defined as 
\[ d_{bL}(\mu,\nu) := \sup \big\{ \int f d\mu - \int f d\nu : f \in \mathcal{F}\big\}, \ \mu,\nu \in \mathcal{P}(\RR),\]
where $\mathcal{F}$ is the class of continuous functions $f : \RR \to \RR$ with Lipschitz constant at most $1$ and uniform bound $1$. It is known that $d_{BL}$ is compatible with the weak topology and that as a consequence of Prokhorov's theorem, $(\mathcal{P}(\RR),d_{bL})$ is a complete metric space. Since $d_{bL}\leq \mathscr{W}_1\leq \mathscr{W}_2$,  Lemma \ref{HW} entails that $W\mapsto \upsilon_W$ is uniformly continuous from $(\mathcal{W}, \| \ \|_{L^1})$ to $(\mathcal{P}(\RR), d_{bL})$. Now, let $W\in \mathcal{W}$.  Since $\| W-W^{(C)}\|_{L^1} \to 0$ as $n\to+\infty$, it follows that $(W^{(C)})_{C>0}$ is a Cauchy sequence for the $L^1$-norm and therefore by uniform continuity $(\upsilon_{W^{(C)}})_{C>0}$ is also a Cauchy sequence in $(\mathcal{P}(\RR),d_{bL})$. As $(\mathcal{P}(\RR),d_{bL})$ is a complete metric space, $(\upsilon_{W^{(C)}})_{C>0}$ is indeed convergent for the weak topology. This shows that $W\mapsto \upsilon_W$ is well-defined on $\mathcal{W}$. 

Now, to show the continuity of the map $W\mapsto \upsilon_W$ for the topology induced by cut norm, let $(U_n)_{n\in \NN}$ be a sequence of $\mathcal{W}$ converging to $U\in\mathcal{W}$ in cut norm. Define for any $C>0$ the kernels $U_{n,C}$ and $\widehat{U}_{n,C}$ as
\[ \widehat{U}_{n,C}(x,y) = U_n(x,y) \Car_{d_U(x)\leq C} \Car_{d_U(y)\leq C},  ,\]
\[ {U}_{n,C}(x,y) = \widehat{U}_{n,C}(x,y) \Car_{d_{U_n}(x)\leq 2C} \Car_{d_{U_n}(y)\leq 2C},  (x,y) \in [0,1]^2.\]
In order to leverage the continuity of $W\mapsto \upsilon_{W}$ on $\mathcal{X}_C$, we will first show that for any $C>0$ and $n\in\NN$,
\begin{equation} \label{claimcont} \|U_{n,C} - U^{(C)}\|_\Box \leq 5\|U_n-U\|_\Box, \end{equation}
\begin{equation} \label{claimcont2} \|U_n-U_{n,C}\|_{L^1} \leq 8\|U_n-U\|_\Box + 4\int d_U \Car_{d_U \geq C} d\lambda.\end{equation}
To prove \eqref{claimcont}, we first note that by definition of the cut norm, we have $\|\widehat{U}_{n,C} - U^{(C)}\|_\Box \leq \|U_n-U\|_\Box$. Thus, it suffices to show that $\|U_{n,C}-\widehat{U}_{n,C}\|_\Box \leq 4\|U_n-U\|_\Box$. Define the event $E_n:= \{d_U\leq C, d_{U_n}>2C\}$. Using the bound $|\widehat{U}_{n,C}(x,y)-U_{n,C}(x,y)| \leq U_n(x,y)[\Car_{E_n}(x) + \Car_{E_n}(y)]$ for any $(x,y) \in [0,1]^2$, we obtain that 
\begin{equation}\label{interm} \|\widehat{U}_{n,C}-U_{n,C}\|_{L^1} \leq 2\int  d_{U_n} \Car_{E_n} d\lambda.\end{equation}
Now, on $E_n$, we have $d_{U_n}\leq 2(d_{U_n}-d_U)$. Therefore, using \eqref{interm}, we get that $\|\widehat{U}_{n,C}-U_{n,C}\|_{L^1}\leq 4\|U_n-U\|_\Box$, which ends the proof of the claim \eqref{claimcont} since $\| \ \|_\Box \leq \| \ \|_{L^1}$. Moving on to showing \eqref{claimcont2}, we note that $|U_n(x,y)-U_{n,C}(x,y)|\leq U_n(x,y)[ \Car_{d_U(x)>C}+\Car_{d_U(y)>C} +\Car_{d_{U_n}(x)>2C} +\Car_{d_{U_n}(y) >2C}]$ for any $(x,y) \in [0,1]^2$. Thus, 
\begin{align} \|U_n-U_{n,C}\|_{L^1} &\leq 2\int d_{U_n} \Car_{\{d_{U_n} >2C\}} d\lambda+2\int d_{U_n} \Car_{\{d_{U} >C\}} d\lambda \nonumber \\
&\leq 2\int d_{U_n} \Car_{E_n}d\lambda +4 \int d_{U_n} \Car_{\{d_{U_n}>2C, d_U>C\}} d\lambda. \label{claim2norm}
\end{align}
On the one hand, as we saw earlier, $d_{U_n} \leq 2(d_{U_n} -d_U)$ on $E_n$ so that $\int d_{U_n} \Car_{E_n} d\lambda \leq 2\|U_n-U\|_\Box$. On the other hand, 
\[ \int d_{U_n} \Car_{\{d_{U_n}>2C, d_U>C\}} d\lambda \leq \|U_n-U\|_\Box + \int d_U \Car_{d_U>C} d\lambda.\]
Plugging this inequality in \eqref{claim2norm}, this yields \eqref{claimcont2}. Now, using the triangle inequality, we write 
\begin{equation}\label{triangle} d(\upsilon_{U_n},\upsilon_{U}) \leq d(\upsilon_{U_n},\upsilon_{U_{n,C}}) + d(\upsilon_{U_{n,C}},\upsilon_{U^{(C)}}) + d(\upsilon_{U^{(C)}},\upsilon_U),\end{equation}
for any $n\in \NN$ and $C>0$. 
On the one hand,  by \eqref{claimcont2} and Lemma \ref{HW} we deduce for any $C>0$ that $\limsup_{n\to+\infty} \mathscr{W}_2(\upsilon_{U_n},\upsilon_{U_{n,C}}) \leq (4\int d_U \Car_{d_U\geq C} d\lambda)^{1/2}$. On the other hand, \eqref{claimcont} shows that $(U_{n,C})_{n\in \NN}$ converges to $U^{(C)}$ in $(\mathcal{X}_C,d_\Box)$ for any $C>0$. By Proposition \ref{conti}, it follows that $(\upsilon_{U_{n,C}})_{n\in \NN}$ converges weakly to $\upsilon_{U^{(C)}}$ for any $C>0$.  Therefore, letting $n\to +\infty$ in \eqref{triangle} and using that $d\leq \mathscr{W}_2$ by \eqref{ineqd}, it yields 
\[ \limsup_{n\to+\infty} d(\upsilon_{U_n},\upsilon_U) \leq (4\int d_U \Car_{d_U\geq C} d\lambda)^{1/2} + d(\upsilon_{U^{(C)}},\upsilon_U).\]
Letting finally $C\to+\infty$ entails that $(\upsilon_{U_n})_{n\in\NN}$ converges weakly to $\upsilon_U$.

\end{proof}


Finally, using the fact that the map $W\mapsto \upsilon_W$ is invariant by the action of ${S}_{[0,1]}$ and  that it is continuous for $d_\Box$, we deduce in the following lemma that it defines a continuous map on the quotient space $\widetilde{\mathcal{W}}$ endowed with the cut metric.

\begin{Pro}\label{conttilde}
Let for any $\widetilde{W}\in \widetilde{\mathcal{W}}$, $\upsilon_{\widetilde{W}}:= \upsilon_W$ where $W\in \widetilde{W}$. The map $\widetilde{W}\mapsto \upsilon_{\widetilde{W}}$ is well-defined on $\widetilde{\mathcal{W}}$ and continuous for the cut metric.

\end{Pro} 
\begin{proof}To prove that $\widetilde{W} \in \widetilde{\mathcal{W}} \mapsto \upsilon_{\widetilde{W}}$ is well-defined amounts to show that whenever $W,W'\in\mathcal{W}$ and $\delta_\Box(W,W') = 0$, one has $\upsilon_W=\upsilon_{W'}$. 
We first  show that the QVE measure of a kernel is invariant by relabelling, in the sense that for any $\phi \in {S}_{[0,1]}$ and  $W\in \mathcal{W}$, we have 
\begin{equation} \label{invQVE} \upsilon_{W^\phi}=\upsilon_{W}.\end{equation} Since $(W^\phi)^{(C)} = (W^{(C)})^\phi$ for any $\phi \in {S}_{[0,1]}$, $C>0$, and $W\in \mathcal{W}$, it is sufficient, by definition of the QVE measure of a kernel, to show that $\upsilon_{W^\phi} = \upsilon_W$ for any $W\in \mathcal{X}$ and $\phi \in {S}_{[0,1]}$. Now, let $\phi \in {S}_{[0,1]}$ and $W\in \mathcal{X}$.  If for any $z\in \mathbb{H}$, $m(z):=(m(z,x))_{x\in[0,1]}$, respectively $\widetilde{m}(z):=(\widetilde{m}(z,x))_{x\in [0,1]}$, is the unique solution in $\mathcal{B}^+$ to the QVE \eqref{QVE} associated to $W$, respectively $W^\phi$, then one can check that $(\widetilde{m}(z,\phi^{-1}(x)))_{x\in [0,1]}$ solves the QVE associated to $W$. By unicity, this implies that $\widetilde{m}(z,\phi^{-1}(x)) = m(z,x)$ for any $x\in [0,1]$, $z\in \mathbb{H}$. As a result, if $g$, respectively $\widetilde{g}$, is the Stieltjes transform of $\upsilon_W$, respectively $\upsilon_{W^\phi}$, then $g(z) = \int_0^1 m(z,x) dx = \int_0^1 \widetilde{m}(z,\phi^{-1}(x)) dx = \int_0^1\widetilde{m}(z,x) dx = \widetilde{g}(z)$ for any $z\in \mathbb{H}$. This implies that $\upsilon_W=\upsilon_{W^\phi}$.

Now, if  $W,W'\in \mathcal{W}$ are such that $\delta_\Box(W,W')=0$, then it means that there exists $(\phi_n)_{n\in \NN}$ a sequence in ${S}_{[0,1]}$ such that $({W'}^{\phi_n})_{n\in \NN}$ converges in cut norm to $W$. By Proposition \ref{defcont} and using \eqref{invQVE}, we have $\upsilon_W=\lim_{n\to+\infty} \upsilon_{{W'}^{\phi_n}} = \upsilon_{W'}$. This shows that the map $\widetilde{W} \in \widetilde{\mathcal{W}} \mapsto \upsilon_{\widetilde{W}}$ is well-defined.
Now, since $W\mapsto \upsilon_W$ is continuous on $(\mathcal{W},d_\Box)$ and $\upsilon_W=\upsilon_{W'}$ whenever $W,W' \in \mathcal{W}$ are such that $\delta_\Box(W,W')=0$, it follows that $\widetilde{W}\mapsto \upsilon_{\widetilde{W}}$ is continuous for the quotient topology on $\widetilde{\mathcal{W}}$, that is, with respect to $\delta_\Box$. 
\end{proof}

\section{Proofs of Theorem \ref{main} and Lemma \ref{zeroI}}

\subsection{Proof of Theorem \ref{main}}
By Proposition \ref{upper} we know that the sequence $(\widetilde{W}_n)_{n\in\NN}$ satisfies a LDP in $(\widetilde{\mathcal{W}},\delta_\Box)$ with speed $n^2p$ and good rate function $H$ defined in \eqref{defHW}. Since the map $\widetilde{W} \mapsto \upsilon_{\widetilde{W}}$ is continuous from $(\widetilde{\mathcal{W}},\delta_\Box)$ to $\mathcal{P}(\RR)$ endowed with the weak topology  by Proposition \ref{conttilde},  it follows from the contraction principle (see \cite[Theorem 4.2.1]{DZ}) that $(\upsilon_{\widetilde{W}_n})_{n\in \NN}$ satisfies a LDP with speed $n^2p$ and rate function $I_L$ defined as 
\begin{equation} \label{defItilde} I_L(\mu) :=\inf \{ H(\widetilde{W}) : \upsilon_{\widetilde{W}} =\mu, \widetilde{W}\in \widetilde{\mathcal{W}}\}, \ \mu \in \mathcal{P}(\RR).\end{equation}
Since for any $\widetilde{W}\in \widetilde{\mathcal{W}}$ and $W'\in \widetilde{W}$, $H(\widetilde{W})=H(W')$ and $\upsilon_{\widetilde{W}}= \upsilon_{W'}$ by definition, we deduce that $I_L$ is indeed represented by the variational problem \eqref{defJ}. Moreover, as $H$ is a good rate function on $(\widetilde{\mathcal{W}},\delta_\Box)$ and $\widetilde{W}\mapsto \upsilon_{\widetilde{W}}$ is continuous, the infimum  \eqref{defItilde} defining $I_L(\mu)$ is achieved for any $\mu$ such that $I_L(\mu)<+\infty$.

\subsection{Proof of Lemma \ref{zeroI}}
Let $\mu_{sc}$ denote the semicircle law. It is well-known that its Stieltjes transform $m$ satisfies the functional equality $-\frac{1}{m(z)} = z+ m(z)$ for any $z\in \mathbb{H}$ (see for example \cite[(2.12)]{BGK}). Thus,  $\mu_{sc}$ is precisely the QVE measure of the kernel constant equal to $1$. Since $h_L(1)=0$ by Lemma \ref{prophL}, we deduce that $I_L(\mu_{sc}) =0$. Assume now that $I_L(\mu) =0$ for some $\mu \in \mathcal{P}(\RR)$. By Theorem \ref{main} we know that the infimum defining $I_L$ is achieved, say at $W\in \mathcal{W}$, meaning that $\upsilon_W =\mu$ and $H(W)=0$. This entails that $h_L\circ W =0$ a.s. Since $h_L$ vanishes only at $1$ by Lemma \ref{prophL}, it follows that $W=1$ a.s. As a result $\mu = \mu_{sc}$.

\appendix
\section{} \label{appendix}

\begin{proof}[Proof of Lemma \ref{rkineq}]
We first note that for any $E$ Borel subset of $[0,1]$ one can find $\phi \in {S}_{[0,1]}$ such that $\phi^{-1}(E)= [0,\lambda(E)]$ (see for example \cite[Theorem A.7]{Janson}). Since $\upsilon_{W^\phi}=\upsilon_W$ for any $W\in\mathcal{X}$ by \eqref{invQVE},  it suffices to prove the statement in the case where $E$ is an interval of the type $[0,t]$ for some $t\in (0,1)$. 
 In a second step we make a further reduction  of the statement to the case where both $W$ and $W'$ are stepped kernels. To this end, fix $W,W'\in \mathcal{X}$ and $t \in (0,1)$ such that $W(x,y) = W'(x,y)$ for almost $(x,y)$ in $[t,1]\times [t,1]$. Set for any $k\in \NN$, $k\geq 1$, $U_k: = W_{\mathcal{P}_k}+1/k$  and $U'_k:=W'_{\mathcal{P}_k}+1/k$, where $W_{\mathcal{P}_k}$ and $W'_{\mathcal{P}_k}$ are the stepped kernels associated to the partition $\mathcal{P}_k:=\{(\frac{i-1}{k},\frac{i}{k}] : i=1,\ldots,k\}$ defined in \eqref{defstepped}  (Adding the small perturbation $k^{-1}$ allows the kernels $U_k$ and $U_k'$ to be bounded away from zero, a fact that we will use later). We will show that for any $k\in \NN$, $k\geq 1$, 
\begin{equation} \label{claimCauchy} d_{\text{KS}}(\upsilon_{U_k},\upsilon_{U_k'}) \leq 2t + \frac{4}{k}.\end{equation}
Assume for the moment that the above inequality is true. As $(U_k)_{k\in \NN}$ and $(U'_k)_{k\in \NN}$ converge respectively to $W$ and $W'$ in $L^1$ norm, therefore as well for the cut norm, we deduce using Proposition \ref{conti} that $\upsilon_{U_k}$ and $\upsilon_{U_k'}$ converge weakly respectively to $\upsilon_W$ and $\upsilon_{W'}$. Thus, taking the limit as $k\to+\infty$ in \eqref{claimCauchy} and using that $d\leq d_{KS}$ by \eqref{ineqd}, we obtain the claimed result.

We now move on to prove \eqref{claimCauchy}. Denote by $S$ and $S'$ the $k\times k$ symmetric matrices associated to the stepped kernels $U_k$ and $U_k'$, that is, let $S_{\ell m}$ and $S'_{\ell m}$ be the values of $U_k$ and $U_k'$ respectively on $(\frac{\ell-1}{k},\frac{\ell}{k}]\times (\frac{m-1}{k},\frac{m}{k}]$ for any $1\leq \ell,m \leq k$. Consider $\Gamma$ to be a $nk\times nk$ symmetric matrix such that $\{\Gamma_{ij}\}_{i\leq j}$ are i.i.d. Rademacher random variables and define the $nk \times nk$ random matrices $Y_n$ and $Y'_n$ by 
\begin{equation} \label{cvmesemp} {Y_n}(i,j) =  \sqrt{\frac{S_{\lceil i/n \rceil, \lceil j/n \rceil}}{nk}} \Gamma_{ij}, \ Y'_n(i,j) =  \sqrt{\frac{S'_{\lceil i/n \rceil, \lceil j/n \rceil}}{nk}} \Gamma_{ij}, \ 1\leq i,j\leq nk. \end{equation}
Since $W,W'\in \mathcal{X}$, we have that $U_k,U_k' \in \mathcal{X}$ as well. It follows that \[ \sup_{n\geq 1} \sup_{1\leq i \leq nk} \frac{1}{nk}\sum_{j=1}^{nk} S_{\lceil i/n \rceil, \lceil j/n \rceil} = \sup_{1\leq \ell \leq k} \frac{1}{k} \sum_{m=1}^k S_{\ell m} = \|d_{U_k}\|_{L^\infty}< +\infty, \]
\[  \sup_{n\geq 1} \sup_{1\leq i \leq nk} \frac{1}{nk}\sum_{j=1}^{nk} S'_{\lceil i/n \rceil, \lceil j/n \rceil} =\sup_{1\leq \ell \leq k} \frac{1}{k} \sum_{m=1}^k S'_{\ell m} =\|d_{U_k'}\|_{L^\infty} <+\infty,\]
As by construction, $U_k(x,y) \geq 1/k$ and $U_k'(x,y) \geq 1/k$ for any $(x,y) \in [0,1]^2$, we also have that $ \inf_{1\leq \ell,m\leq k} S_{\ell m} \geq 1/k$ and $\inf_{1\leq \ell,m\leq k} S'_{\ell m}\geq 1/k$.
By \cite[Theorem 1.1]{Girko}, we obtain that almost surely
\begin{equation} \label{cvmes1} d_{\text{KS}}(\mu_{Y_n},\upsilon_{U_k}) \underset{n\to+\infty}{\longrightarrow} 0, \ d_{\text{KS}}(\mu_{Y'_n},\upsilon_{U_k'}) \underset{n\to+\infty}{\longrightarrow} 0.\end{equation}
We now exploit our assumption on the kernels $W$ and $W'$. Fix some $\ell_0\geq 2$ such that $\frac{\ell_0-2}{k} < t \leq  \frac{\ell_0-1}{k}$. Note that since $W$ and  $W'$ coincides almost surely on $[t,1]\times [t,1]$, we have that $U_k(x,y) = U_k'(x,y)$ whenever $x>\frac{\ell_0-1}{k}$ and $y>\frac{\ell_0-1}{k}$. Thus, $S_{\ell m} =S'_{\ell m}$ for any $\ell,m\geq \ell_0$, and as a consequence $Y_{ij} = Y'_{ij}$ for any $i,j \geq \ell_0 n$. In particular, $Y-Y'$ has rank at most $2\ell_0 n$. By Cauchy interlacing inequality (see \cite[Theorem 1.43]{BS}), we have
\begin{equation} \label{Cauchyfinite} d_{\text{KS}}(\mu_{Y_n},\mu_{Y'_n}) \leq \frac{2\ell_0}{k}\leq  2t+\frac{4}{k}.\end{equation}
Taking the limit as $n$ goes to $+\infty$, we obtain using \eqref{cvmes1} the inequality \eqref{claimCauchy}. 

\end{proof}

\begin{proof}[Proof of Lemma \ref{HW}] We first reduce the statement to the case of stepped kernels. Let $W,W'\in \mathcal{X}$ and define $U_k$, $U_k'$ for any $k\in \NN$ as in the proof of Lemma \ref{rkineq}. We will show that for any $k\in \NN$, $k\geq 1$, 
\begin{equation} \label{WHstepped} \mathscr{W}_2(\upsilon_{U_k},\upsilon_{U'_k}) \leq  \|U_k-U'_k\|_{L^1}^{1/2}.\end{equation}
Assume for the moment that \eqref{WHstepped} holds. Since $(U_k)_{k\geq 1}$ and $(U_k')_{k\geq 1}$ converge for the  $L^1$ norm to respectively $W$ and $W'$, therefore as well for the cut norm, it follows by Proposition \ref{conti} that $(\upsilon_{U_k})_{k\geq 1}$ and $(\upsilon_{U_k'})_{k\geq 1}$ converge weakly respectively to $\upsilon_W$ and $\upsilon_{W'}$. In addition, by Lemma \ref{momentuW} we have for any $k\geq 1$, 
\[ \int \tau^2 d\upsilon_{U_k}(\tau) = \int_{[0,1]^2} U_k(x,y) dx dy \underset{k\to+\infty}{\longrightarrow} \int_{[0,1]^2} W(x,y)dx dy = \int \tau^2 d\upsilon_W(\tau), \]
where we used again the fact that $(U_k)_{k\geq 1}$ converges in $L^1$ norm to $W$. By \cite[Theorem 6.9]{Villani}, we deduce that $(\upsilon_{U_k})_{k\geq 1}$ converges for the $L^2$ Wasserstein metric to $\upsilon_W$. Similarly, we have that $\mathscr{W}_2(\upsilon_{U_k'},\upsilon_{W'}) \to 0$ as $k\to +\infty$. Thus, taking the limit in \eqref{WHstepped} when $k\to +\infty$ gives the claimed inequality. 

It remains to prove \eqref{WHstepped}. As in the proof of Lemma \ref{rkineq}, we define $Y_n$ and $Y_n'$ $nk\times nk$ symmetric random matrices as in \eqref{cvmesemp}. We know from \eqref{cvmes1} that almost surely $(\mu_{Y_n})_{n\in \NN}$ and $(\mu_{Y_n'})_{n\in \NN}$ converge weakly respectively to $\upsilon_{U_k}$ and $\upsilon_{U_k'}$. Now, observe that for any $n\geq 1$,
\begin{align*}
 \int \tau^2 d \mu_{Y_n}(\tau) &= \frac{1}{nk} \sum_{1\leq i,j\leq nk} Y_n(i,j)^2 = \frac{1}{(nk)^2} \sum_{1\leq i,j\leq nk} S_{\lceil i/n \rceil, \lceil j/n\rceil} \\
& = \int_{[0,1]^2} U_k(x,y) dx dy.
\end{align*}
Since $\int \tau^2 d\upsilon_{U_k}(\tau) = \int_{[0,1]^2} U_k(x,y) dx dy$ by Lemma \ref{momentuW}, it follows again by \cite[Theorem 6.9]{Villani} that $(\mu_{Y_n})_{n\geq 1}$ converges in $L^2$ Wasserstein distance to $\upsilon_{W_k}$ almost surely. Similarly, $\mathscr{W}_2(\mu_{Y_n'},\upsilon_{U_k'}) \to 0$ as $n\to+\infty$. <by Hoeffman-Wielandt inequality (see \cite[Lemma 2.1.19]{AGZ})  we have for any $n\geq 1$, 
\[ \mathscr{W}_2(\mu_{Y_n},\mu_{Y_n'})^2 \leq \frac{1}{kn} {\rm Tr}(Y_n-Y'_n)^2 = \frac{1}{(kn)^2} \sum_{1\leq i,j\leq kn} \big|\sqrt{S_{\lceil \frac{i}{n} \rceil \lceil \frac{j}{n}\rceil}}-\sqrt{S'_{\lceil \frac{i}{n} \rceil \lceil \frac{j}{n}\rceil}}\big|^2 .\]
Using that $|\sqrt{x}-\sqrt{y}|\leq\sqrt{|x-y|}$ for any $x,y \geq 0$, we find that
\begin{align*}
 \mathscr{W}_2(\mu_{Y_n},\mu_{Y_n'})^2& \leq  \frac{1}{(kn)^2} \sum_{1\leq i,j\leq kn} \big|S_{\lceil \frac{i}{n} \rceil \lceil \frac{j}{n}\rceil}-S'_{\lceil \frac{i}{n} \rceil \lceil \frac{j}{n}\rceil}\big|= \frac{1}{n^2} \sum_{1\leq \ell,m \leq k} |S_{\ell m} - S'_{\ell m}|\\
& = \|U_k-U_k'\|_{L^1}.
\end{align*}
Taking the limit in the above inequality as $n\to +\infty$, gives then the claim \eqref{WHstepped}.

\end{proof}

\begin{Lem}\label{chgmes}
Let $\PP$ and $\QQ$ be two probability measures defined on the same $\sigma$-algebra $\mathcal{A}$. 
For any $E\in \mathcal{A}$ such that $\QQ(E)>0$, 
\[ \PP( E) \geq  \QQ(E)\exp\Big( - \frac{1}{\QQ(E)} \big( H(\QQ\mid \PP) +e^{-1}\big) \Big),\]
where $H(\QQ\mid \PP)$ is defined in \eqref{defentropy}
\end{Lem}

\begin{proof}Without loss of generality, we can assume that $\QQ$ is absolutely continuous with respect to $\PP$, otherwise the inequality is trivially satisfied.
Using Jensen's inequality, we get 
\[ \PP(E) = \QQ(E)\EE_\QQ\big( \frac{\Car_E}{\QQ(E)} e^{-\log \frac{d\QQ}{d\PP}}\big) )\geq \QQ(E) \exp\Big(- \frac{1}{\QQ(E)} \EE_\QQ\big( \Car_{E} \log \frac{d\QQ}{d\PP}\big) \Big) .\]
As $x \log x \geq -e^{-1}$ for any $x>0$, it follows that 
\[ \EE\big( \Car_{E} \frac{d\QQ}{d\PP} \log \frac{d\QQ}{d\PP}\big)  \leq \EE\big( \frac{d\QQ}{d\PP} \log \frac{d\QQ}{d\PP}\big)  +e^{-1},\]
which ends the proof.
\end{proof}

\bibliographystyle{plain}
\bibliography{main.bib}{}

\end{document}